\documentclass[10pt]{article}
\usepackage{lmodern}
\usepackage{amsmath}
\usepackage[T1]{fontenc}
\usepackage[utf8]{inputenc}
\usepackage{authblk}
\usepackage{amsfonts}
\usepackage{graphicx}
\usepackage{amssymb}
\usepackage[english]{babel}
\usepackage{color}
\usepackage{amsthm}
\usepackage{graphicx}
\usepackage{mathrsfs}
\usepackage{makecell}
\usepackage{microtype}
\usepackage{mathscinet}
\usepackage{array}
\usepackage{multirow}
\usepackage{enumerate}
\usepackage[cal=boondoxo,bb=ams]{mathalfa}
\usepackage{hyperref}
\hypersetup{hidelinks}
\usepackage{booktabs}
\usepackage{arydshln}

\newtheorem{theorem}{Theorem}[section]
\newtheorem{prop}{Proposition}[section]
\newtheorem{lemma}{Lemma}[section]

\newtheorem{remark}{Remark}[section]

\newcommand{\ml}{\mathcal}
\newcommand{\mb}{\mathbb}

\DeclareMathOperator{\lin}{lin}
\DeclareMathOperator{\intt}{int}
\DeclareMathOperator{\extt}{ext}
\DeclareMathOperator{\bdd}{bdd}

\title{On the Cauchy problem for acoustic waves in hereditary fluids: decay properties and inviscid limits}
\author[1]{Wenhui Chen\thanks{Wenhui Chen (wenhui.chen.math@gmail.com)}}
\affil[1]{School of Mathematics and Information Science, Guangzhou University, 510006 Guangzhou, China}
\date{}

\begin{document}

		\maketitle
		\begin{abstract}
			\medskip
	This manuscript considers the viscous/inviscid Moore-Gibson-Thompson (MGT) equations with memory of type I in the whole space $\mb{R}^n$. For one thing, associating with a new condition on initial data, we derive the optimal $L^2$ estimates and the optimal leading term of the acoustic velocity potential for large time, in which we analyze different contributions from viscous, thermally relaxing, as well as hereditary fluids on large time asymptotic behavior for the acoustic waves models. For another, via the multi-scale analysis and energy methods in the Fourier space, we demonstrate the $L^{\infty}$ inviscid limits in the sense of the diffusivity of sound tending to zero, which match our WKB expansion of the solution. Finally, we give a further application of our results on large time behavior for the nonlinear Jordan-MGT equation in viscous hereditary fluids.
			\\
			
			
			\noindent\textbf{AMS Classification (2020)}  35L30, 74D05, 35B40, 35B20
		\end{abstract}
\fontsize{12}{15}
\selectfont

\section{Introduction}\label{Section-Introduction}\setcounter{equation}{0}
$\ \ \ \ $In the present paper, we study the Cauchy problem for the viscous/inviscid Moore-Gibson-Thompson (MGT) equations with memory of type I, arising from the linearized models for acoustic waves propagation in viscous/inviscid thermally relaxing hereditary fluids, namely,
\begin{align}\label{MGT-Memory}
	\begin{cases}
		\tau\psi_{ttt}+\psi_{tt}-\Delta\psi-(\delta+\tau)\Delta\psi_t+g^{\tau}\ast\Delta\psi=0,&x\in\mb{R}^n,\ t>0,\\
		(\psi,\psi_t,\psi_{tt})(0,x)=(\psi_0,\psi_1,\psi_2)(x),&x\in\mb{R}^n,
	\end{cases}
\end{align}
with the thermal relaxation $\tau>0$ and the diffusivity of sound $\delta\geqslant 0$, where the unknown function $\psi=\psi(t,x)\in\mb{R}$ is referred to the acoustic velocity potential in a given flow with memory effect. The memory term of type I is expressed by
\begin{align*}
g^{\tau}\ast\Delta\psi:=\int_0^tg^{\tau}(t-\eta)\Delta\psi(\eta,x)\mathrm{d}\eta
\end{align*}
with the memory kernel $g^{\tau}=g^{\tau}(t)$.  In the thermally relaxing media of our consideration, e.g. the high-frequency waves in liquids and gases, the memory kernel is typically represented by the exponential decay form (cf. \cite[Chapter 1]{Naugolnykh-Ostrovsky=2000} and \cite[Section 1]{Dell-Pata=2017}) as follows:
\begin{align}\label{Memory-Kernel}
g^{\tau}(t):=m\mathrm{e}^{-\frac{t}{\tau}}\ \ \mbox{with}\ \ m\tau<1,
\end{align}
carrying the small relaxation parameter $m>0$. According to the physical nature of this small parameter $m$, which is related to the sound velocity in high-frequency limit, the last condition $m\tau<1$ always holds. Our main purpose is to deeply understand asymptotic behavior for the Cauchy problem \eqref{MGT-Memory} with the memory kernel \eqref{Memory-Kernel} in the viscous case $\delta>0$ and the inviscid case $\delta=0$. Particularly, via investigating large time behavior and inviscid limits (in the sense of the diffusivity of sound $\delta$ tending to zero), we are interested in different contributions from viscous, thermally relaxing, as well as hereditary fluids on the acoustic waves models.

Since the beginning of this century, theoretical studies of acoustic waves are widely applied in medical and industrial uses of high-intensity
	ultra sound, e.g. medical imaging and therapy, ultrasound cleaning and welding (cf. \cite{Abramov-1999,Dreyer-Krauss-Bauer-Ried-2000,Kaltenbacher-Landes-Hoffelner-Simkovics-2002} and references given therein). In order to describe propagation of sound in viscous thermally relaxing fluids, the classical researches of nonlinear acoustic waves with second sound phenomenon (cf. the pioneering works \cite{Blackstock-1963,Hamilton-Blackstock-1998,Jordan-2014}) deduced some approximated models of the fully compressible Navier-Stokes-Cattaneo system with the pressure-density state equation under irrotational flows, in which the Maxwell-Cattaneo law of heat conduction is to amend the unphysical phenomenon of an infinite signal speed. Among them, via the Lighthill scheme of approximations to retain the first order and second order terms in the sense of small perturbations around the constant equilibrium state, the following well-known viscous/inviscid Jordan-MGT equations (see \cite{Kaltenbacher-Lasiecka-Pos-2012,Kaltenbacher-Nikolic-2019,Racke-Said-2020,B-L-2020,Kaltenbacher-Niko-2021,Chen-Takeda=2023} and references therein) arise:
\begin{align}\label{JMGT}
	\tau\varphi_{ttt}+\varphi_{tt}-\Delta\varphi-(\delta+\tau )\Delta\varphi_t=\partial_t\left(\frac{B}{2A}|\varphi_t|^2+|\nabla\varphi|^2\right)
\end{align}
with $\tau>0$ and $\delta\geqslant0$,
carrying the coefficients of nonlinearity satisfying $B/A>0$. Remark that $\delta>0$ is always referred to the viscous case due to its viscous dissipation originating from the Navier-Stokes equations, analogously, $\delta=0$ may stand for the inviscid case. The corresponding linearizations of the Jordan-MGT equation \eqref{JMGT}, i.e. the MGT equations (named after F.K. Moore, W.E. Gibson \cite{MooreGibson1960} in 1960 and P.A. Thompson \cite{Thompson1972} in 1972), are addressed by 
\begin{align*}
	\tau\psi_{ttt}+\psi_{tt}-\Delta\psi-(\delta+\tau )\Delta\psi_t=0.
\end{align*}
Actually, this model equation firstly appeared in a very old paper of G.G. Stokes \cite{Stokes=1851} in 1851. Some qualitative properties for the MGT equation in bounded or unbounded domains are deeply investigated by \cite{Kaltenbacher-Lasiecka-Marchand-2011,Marchand-McDevitt-Triggiani-2012,Conejero-Lizama-Rodenas-2015,Pellicer-Said-Houari=2019,Chen-Ikehata=2021,Chen-Takeda=2023,Chen=2023} and references therein, including well-posedness, chaotic dynamics, propagation of singularities, large time behavior (optimal estimates and asymptotic profiles), inviscid limits and singular limits associated with initial layers.

As relaxation processes appear in high-frequency waves, the acoustic pressure may depend on the medium density at all prior time. Particularly, these relaxation processes can occur when there are inhomogeneities in the propagation region, for instance, through excitation of molecular degrees of freedom or impurity effect in the fluids (cf. \cite[Chapter 1]{Naugolnykh-Ostrovsky=2000}). In such cases, the pressure-density state equation is not satisfied exactly but up to a term that involves the history of the process, namely, with a memory term. By using the Lighthill scheme of approximations again, the next Jordan-MGT equations with memory of type I occurs:
\begin{align}\label{JMGT-Memory}
	\tau\varphi_{ttt}+\varphi_{tt}-\Delta\varphi-(\delta+\tau )\Delta\varphi_t+g^{\tau}\ast\Delta\varphi=\partial_t\left(\frac{B}{2A}|\varphi_t|^2+|\nabla\varphi|^2\right)
\end{align}
with $\tau>0$ and $\delta\geqslant0$, whose well-posedness and decay properties have been studied by \cite{Caixeta-Lasiecka-Domingos-Valeria=2016,Lasiecka=2017,Lizama-Zamorano=2019,Nikolic-Said=2021,Nikolic-Said=2021-02,Nikolic-Said=2021-03,Boul-Chou-Ouchen=2021}. Concerning its linearized models, namely,
\begin{align}\label{MGT-Memory-Intro}
\tau\psi_{ttt}+\psi_{tt}-\Delta\psi-(\delta+\tau)\Delta\psi_t+g^{\tau}\ast\Delta\psi=0,
\end{align}
the authors of \cite{Lasiecka-Wang-2016,Lasiecka-Wang=2015} studied decay estimates of an energy term in the viscous case $\delta>0$ by assuming some decay properties for the memory kernel $g^{\tau}$, moreover, the authors of \cite{Dell-Lasiecka=Pata=2016} derived polynomial stabilities in the inviscid case $\delta=0$. Later, \cite{Alves-Caixeta-Silva-Jorge=2018} proved the uniform stability of the viscous model \eqref{MGT-Memory-Intro} by the linear semigroup theory. For the Cauchy problem \eqref{MGT-Memory} with some assumptions on the memory kernel, the authors of \cite{Bou-Said=2021} derived decay estimates for the energy term as follows:
\begin{align}\label{Bou-Said}
&\|(\psi_t+\tau \psi_{tt},\nabla\psi+\tau \nabla\psi_t,\nabla \psi_t)\|_{L^2}\notag\\
&\lesssim\begin{cases}
(1+t)^{-\frac{n}{4}}\|(\psi_1+\tau \psi_{2},\nabla\psi_0+\tau \nabla\psi_1,\nabla \psi_1)\|_{(L^2\cap L^1)^{1+2n}}&\mbox{when}\ \ \delta>0,\\
(1+t)^{-\min\{\frac{n}{4},\frac{\ell}{2}\}}\|(\psi_1+\tau \psi_{2},\nabla\psi_0+\tau \nabla\psi_1,\nabla \psi_1)\|_{(H^{\ell}\cap L^1)^{1+2n}}&\mbox{when}\ \ \delta=0,
\end{cases}
\end{align}
with $\ell\geqslant0$. In other words, the regularity-loss type decay property (i.e. to get decay estimates, we need to lose $\ell$ regularity for initial data) occurs only in the inviscid case $\delta=0$. Their approaches are based on energy
methods in Fourier space associated with Lyapunov functionals.  Concerning different decay hypotheses on the memory kernel, we refer the interested reader to \cite{Liu-Chen=2020,Liu-Chen-Chen=2020,Dell-Lasiecka-Pata=2020,Lach-Messao=2021,Nicaise-Bouna=2021,Zhang=2022,Conti-Liv-Pata=2023} and references given therein. To the best of author's knowledge, asymptotic behavior for the solution itself $\psi$ to the Cauchy problem \eqref{MGT-Memory} is completely unknown. We will partly answer this question via large time behavior and small diffusivity of sound behavior for the acoustic velocity potential $\psi$ in this work.

Our main contributions in the present paper are stated in Section \ref{Section-Results}. To be specific, in Section \ref{Section-MGT-Large-Time}, under the non-trivial condition $P_{\psi_1+\tau\psi_2}\neq0$, by assuming suitable Sobolev regularities for initial data in the Cauchy problem \eqref{MGT-Memory} with $\delta\geqslant 0$, motivated by the reduction procedure to the fourth order (in time) evolution equation, we obtain the optimal growth ($n=1,2$) and decay ($n\geqslant 3$) estimates of the solution $\psi$ for large time $t\gg1$. Especially, we find the new regularity-loss threshold of large time optimal estimates for the inviscid model \eqref{MGT-Memory} with $\delta=0$ which is the dimension $n=6$. Furthermore, with the new condition $Q_{\psi_0,\psi_1,\psi_2}\neq0$ defined in \eqref{Quantity}, via the WKB analysis and the refined Fourier analysis, we derive the optimal leading term in the form of diffusion-waves. These results show different influence of viscous, thermally relaxing, as well as hereditary fluids on large time behavior for the acoustic wave models \eqref{MGT-Memory}, which will be shown at length in Remarks \ref{Rem-Influence-01} and \ref{Rem-Influence-02}. For another, in Section \ref{Section-Inviscid}, by using the multi-scale analysis and energy methods in the Fourier space, we rigorously justify the $L^{\infty}$ inviscid limits, which should be understood in this context as the vanishing sound diffusivity limit (cf. \cite[Section 1]{Kaltenbacher-Niko-2021}), of the suitable energy term globally in-time and of the acoustic velocity potential $\psi$ locally in-time. The convergence rate $\sqrt{\delta}$ of the $L^{\infty}$ inviscid limits exactly matches our WKB expansion of the solution. Finally in Section \ref{Section-Final-Remarks}, we give a further application of our derived results on large time behavior for the viscous Jordan-MGT equation \eqref{JMGT-Memory} in hereditary fluids.

\section{Main results}\label{Section-Results}\setcounter{equation}{0}
\subsection{Result and discussion on the optimal estimates}
$\ \ \ \ $Our first result contributes to the optimal growth ($n=1,2$) and decay ($n\geqslant3$) estimates of the solution $\psi$ to the MGT equations \eqref{MGT-Memory} with $\delta\geqslant 0$ in hereditary fluids with initial data localizing in the function space $\ml{A}_{\delta,\ell}$ as follows:
\begin{align*}
\ml{A}_{\delta,\ell}:=\begin{cases}
	(L^2\cap L^1)\times(L^2\cap L^1)\times (L^2\cap L^1)&\mbox{when}\ \ \delta>0,\\
(H^{\ell}\cap L^1)\times(H^{\max\{\ell-1,0\}}\cap L^1)\times (H^{\max\{\ell-2,0\}}\cap L^1)&\mbox{when}\ \ \delta=0,
\end{cases}
\end{align*}
under suitable regularities $\ell\geqslant0$. That is to say, we need higher Sobolev regularities of initial data in the inviscid case $\delta=0$.
\begin{theorem}\label{Thm-Optimal-Est}
Suppose that initial data $(\psi_0,\psi_1,\psi_2)\in\ml{A}_{\delta,\ell}$ with $\ell>\max\{\frac{n}{2}-1,0\}$ and $\delta\geqslant0$. Then, the solution to the MGT equations \eqref{MGT-Memory} with the memory \eqref{Memory-Kernel} satisfies the following optimal growth/decay estimates:
	\begin{align}\label{Est-Itself}
		\ml{D}_n(t)|P_{\psi_1+\tau\psi_2}|\lesssim\|\psi(t,\cdot)\|_{L^2}\lesssim\ml{D}_n(t)\|(\psi_0,\psi_1,\psi_2)\|_{\ml{A}_{\delta,\ell}}
	\end{align}
for large time $t\gg1$, provided that $P_{\psi_1+\tau\psi_2}\neq0$, in which the time-dependent function $\ml{D}_n(t)$ is defined via
\begin{align}\label{Dnt} 
	\ml{D}_n(t):=\begin{cases}
		\sqrt{t}&\mbox{when}\ \ n=1,\\
		\sqrt{\ln t}&\mbox{when}\ \ n=2,\\
		t^{\frac{1}{2}-\frac{n}{4}}&\mbox{when}\ \ n\geqslant 3.
	\end{cases}
\end{align}
Here, we denote $P_f:=\int_{\mb{R}^n}f(x)\mathrm{d}x$ to be the mean of the function $f=f(x)$.
\end{theorem}
\begin{remark}
The new regularity-loss threshold of large time optimal estimates for the inviscid model \eqref{MGT-Memory} with $\delta=0$ is the dimension $n=6$.  Namely, when $1\leqslant n\leqslant 5$ the solution satisfies the optimal estimates \eqref{Est-Itself} with $L^2\cap L^1$ data; however, when $n\geqslant 6$ it requires $H^{\frac{n}{2}-3+\epsilon_0}\cap L^1$ regularity with a sufficiently small parameter $\epsilon_0>0$ for initial data to estimate the solution $\psi$ optimally in the $L^2$ norm.
\end{remark}
\begin{remark}\label{Rem-Influence-01} We state some influence from the memory term of type I on the MGT equations. Let us firstly recall the Cauchy problem for the viscous/inviscid MGT equations
	\begin{align}\label{MGT}
		\begin{cases}
			\tau\psi_{ttt}+\psi_{tt}-\Delta\psi-(\delta+\tau)\Delta\psi_t=0,&x\in\mb{R}^n,\ t>0,\\
			(\psi,\psi_t,\psi_{tt})(0,x)=(\psi_0,\psi_1,\psi_2)(x),&x\in\mb{R}^n.
		\end{cases}
	\end{align}
	 For one thing, comparing with the corresponding result for the viscous MGT equation \eqref{MGT} carrying $\delta>0$, the recent papers \cite{Chen-Ikehata=2021,Chen-Takeda=2023} derived the same optimal estimates as those in Theorem \ref{Thm-Optimal-Est}. In other words, the additional memory term of type I, i.e. $g^{\tau}\ast\Delta\psi$, with the exponential decay kernel \eqref{Memory-Kernel} does not change the estimates for the acoustic velocity potential. It may be caused by the stronger influence from the viscous term. For another, concerning the inviscid MGT equation \eqref{MGT} with $\delta=0$, later, Theorem \ref{Prop-Inviscid-Optimal-Growth} obtains the optimal growth estimates for the solution in the $L^2$ norm when $n=1,2$. With the aid of additional memory effect, we may realize the optimal decay estimates for higher dimensions, but we require higher Sobolev regularity for initial data. In order to understand this influence of the memory term, we next state the comparison table for the inviscid models.
\begin{table}[h!]
	\begin{center}
		\caption{Influence on the inviscid models with $\delta=0$ from thermally relaxing hereditary fluids}
		\medskip
		\label{tab:table1}
		\begin{tabular}{cccccc} 
			\toprule
			 & Reference & $n=1$ & $n=2$ & $n\geqslant3$ & Regularity of data\\
			\midrule
			Free wave equation &\multirow{2}{*}{\cite{Ikehata=2023,Chen-Takeda=2023}} & \multirow{2}{*}{$\sqrt{t}$} & \multirow{2}{*}{$\sqrt{\ln t}$} & \multirow{2}{*}{--}& \multirow{2}{*}{$L^2\times (L^2\cap L^1)$} \\
			($m=0$, $\tau=0$) & & & & &\\
			\midrule
			Inviscid MGT equation &\multirow{2}{*}{Theorem \ref{Prop-Inviscid-Optimal-Growth}} & \multirow{2}{*}{$\sqrt{t}$} & \multirow{2}{*}{$\sqrt{\ln t}$} & \multirow{2}{*}{--}& \multirow{2}{*}{$L^2\times (L^2\cap L^1)^2$} \\
			($m=0$, $\tau>0$)& & & & & \\
			\midrule
			Inviscid MGT equation &\multirow{3}{*}{Theorem \ref{Thm-Optimal-Est}}  & \multirow{3}{*}{$\sqrt{t}$} & \multirow{3}{*}{$\sqrt{\ln t}$} & \multirow{3}{*}{$t^{\frac{1}{2}-\frac{n}{4}}$}& $(H^{\ell}\cap L^1)\times (H^{\max\{\ell-1,0\}}\cap L^1)$ \\
			with the memory & & & & & $\times(H^{\max\{\ell-2,0\}}\cap L^1)$\\
			($m>0$, $\tau>0$) & & & & & with $\ell>\max\{\frac{n}{2}-1,0\}$\\
			\bottomrule
			\multicolumn{6}{l}{\emph{$*$ All time-dependent coefficients of the $L^2$ estimate in the table are optimal for large time.}}
		\end{tabular}
	\end{center}
\end{table}
\end{remark}

\begin{remark}
Even the energy term $(\psi_t+\tau\psi_{tt},\nabla\psi+\tau\nabla\psi_t,\nabla\psi_t)$ decays polynomially in the $L^2$ norm found by \cite{Bou-Said=2021}, i.e. the decay estimate \eqref{Bou-Said}, the acoustic velocity potential $\psi$ may grow polynomially when $n=1$ or logarithmically when $n=2$.
\end{remark}
\begin{remark}
Concerning the inviscid model \eqref{MGT-Memory} with $\delta=0$, in Theorem \ref{Thm-Optimal-Est} our regularity condition $\ell>\max\{\frac{n}{2}-1,0\}$ is to guarantee the optimal estimates \eqref{Est-Itself}. Actually, the general estimate \eqref{General-Est} holds for any $\ell\geqslant0$. The question of large time optimal estimates for lower regular data with $0\leqslant \ell\leqslant\max\{\frac{n}{2}-1,0\} $ is still open.
\end{remark}
\subsection{Result and discussion on the optimal leading term}
$\ \ \ \ $Next, let us introduce the first large time profile $\psi^{(1,p)}=\psi^{(1,p)}(t,x)$ such that
\begin{align*}
\psi^{(1,p)}(t,x):=J_0(t,x)P_{\psi_1+\tau\psi_2}
\end{align*}
and the second large time profile $\psi^{(2,p)}=\psi^{(2,p)}(t,x)$ such that
\begin{align*}
\psi^{(2,p)}(t,x):=-\nabla J_0(t,x)\circ M_{\psi_1+\tau\psi_2}+J_1(t,x)P_{\psi_0-\tau^2\psi_2}+J_2(t,x)P_{\psi_1+\tau\psi_2},
\end{align*}
where the auxiliary diffusion-waves functions $J_k=J_k(t,x)$ with $k=0,1,2$ are defined via the partial Fourier transform as follows:
\begin{align*}
	J_0&:=\ml{F}^{-1}_{\xi\to x}\left(\frac{\sin(\sqrt{1-m\tau}|\xi|t)}{\sqrt{1-m\tau}|\xi|}\mathrm{e}^{-\frac{1}{2}(\delta+2m\tau^2)|\xi|^2t}\right),\\ J_1&:=\ml{F}^{-1}_{\xi\to x}\left(\cos(\sqrt{1-m\tau}|\xi|t)\mathrm{e}^{-\frac{1}{2}(\delta+2m\tau^2)|\xi|^2t}\right),\\
	J_2&:=\frac{\gamma_0t}{\sqrt{1-m\tau}}\ml{F}^{-1}_{\xi\to x}\left(|\xi|^2\cos(\sqrt{1-m\tau}|\xi|t)\mathrm{e}^{-\frac{1}{2}(\delta+2m\tau^2)|\xi|^2t}\right).
\end{align*}
In the above settings, we denote $M_f:=\int_{\mb{R}^n}xf(x)\mathrm{d}x$ to be the weighted mean of the function $f=f(x)$. Moreover, the symbol $\circ$ stands for the inner product in the Euclidean space. For the sake of simplicity, we take the parameter $\gamma_0=\gamma_0(\tau,\delta,m)$ via
\begin{align}\label{gamma-0}
	\gamma_0:=\frac{(4\tau-\delta)\delta-4m\tau^2(2\delta+4m\tau^2-3)}{8\sqrt{1-m\tau}}.
\end{align}
\begin{remark}\label{Rem-Influence-02}
For the viscous memoryless case $m=0$ with $\delta>0$, i.e. the viscous MGT equation \eqref{MGT}, our first and second large time profiles exactly coincide with those in \cite[Subsection 2.5]{Chen-Takeda=2023}, in which we also correct the misprint of the sign for $\nabla J_0(t,x)$ in the expression of $\psi^{(2,p)}$ with an additional minus. One of the influence from the memory term of type I on large time profiles for the viscous MGT equation \eqref{MGT} with $\delta>0$, actually, is reflected on the quantitative properties. Precisely, there are some additional factors as well as coefficients related to the parameter $m$. The qualitative properties of large time profiles, however, in the viscous model do not be changed. More details in the comparison are stated in the next table for the viscous models.
\begin{table}[h!]
	\begin{center}
		\caption{Influence on the viscous models with $\delta>0$ from thermally relaxing hereditary fluids}
		\medskip
		\label{tab:table2}
		\begin{tabular}{ccc} 
			\toprule
			& Reference & Optimal leading term \\
			\midrule
			Strongly damped wave equation &\multirow{2}{*}{\cite{Ikehata=2014,Ikehata-Ono=2017}} & \multirow{2}{*}{$\ml{F}^{-1}_{\xi\to x}\left(\frac{\sin(|\xi|t)}{|\xi|}\mathrm{e}^{-\frac{1}{2}\delta|\xi|^2t}\right)$} \\
			($m=0$, $\tau=0$)& &\\
			\midrule
			Viscous MGT equation & \multirow{2}{*}{\cite{Chen-Ikehata=2021,Chen-Takeda=2023}} &\multirow{2}{*}{$\ml{F}^{-1}_{\xi\to x}\left(\frac{\sin(|\xi|t)}{|\xi|}\mathrm{e}^{-\frac{1}{2}\delta|\xi|^2t}\right)$}\\
			($m=0$, $\tau>0$)&&\\
			\midrule
			Viscous MGT equation &\multirow{3}{*}{Theorem \ref{Thm-Optimal-Leading}}  & \multirow{3}{*}{$\ml{F}^{-1}_{\xi\to x}\left(\frac{\sin(\sqrt{1-m\tau}|\xi|t)}{\sqrt{1-m\tau}|\xi|}\mathrm{e}^{-\frac{1}{2}(\delta+2m\tau^2)|\xi|^2t}\right)$} \\
			with the memory & &  \\
			($m>0$, $\tau>0$) & & \\
			\bottomrule
			\multicolumn{3}{l}{\emph{$*$ All leading terms in the table fulfill the optimal error estimates with the same decay rate.}}
		\end{tabular}
	\end{center}
\end{table}
\end{remark}
\begin{remark}
We propose a significant influence from the memory term of type I on the inviscid MGT equation \eqref{MGT} with $\delta=0$, whose solution $\psi^{m=0,\delta=0}=\psi^{m=0,\delta=0}(t,x)$ can be explicitly represented in the Fourier space (see \cite[Proof of Proposition 2]{Chen-Palmieri=2020} in detail) via
\begin{align}
\widehat{\psi}^{m=0,\delta=0}&=\left(\frac{\cos(|\xi|t)}{1+\tau^2|\xi|^2}+\frac{\tau|\xi|\sin(|\xi|t)}{1+\tau^2|\xi|^2}+\frac{\tau^2|\xi|^2}{2(1+\tau^2|\xi|^2)}\mathrm{e}^{-\frac{t}{\tau}}\right)\widehat{\psi}_0+\frac{\sin(|\xi|t)}{|\xi|}\widehat{\psi}_1\notag\\
&\quad+\left(\frac{\tau\sin(|\xi|t)}{|\xi|(1+\tau^2|\xi|^2)}-\frac{\tau^2\cos(|\xi|t)}{1+\tau^2|\xi|^2}+\frac{\tau^2}{2(1+\tau^2|\xi|^2)}\mathrm{e}^{-\frac{t}{\tau}}\right)\widehat{\psi}_2.\label{Rep-Inviscid}
\end{align}
That is to say, we may understand large time behavior for the inviscid MGT equation by the free waves structure because of the crucial factors (cf. Appendix \ref{Appendix-A})
\begin{align}\label{W-01}
\cos(|\xi|t)\ \ \mbox{as well as}\ \ \frac{\sin(|\xi|t)}{|\xi|}.
\end{align}
 With the aid of additional memory effect, its large time behavior is dominated by the diffusion-waves structure, that is,
 \begin{align}\label{W-02}
 \cos(\sqrt{1-m\tau}|\xi|t)\mathrm{e}^{-m\tau^2|\xi|^2t}\ \ \mbox{as well as}\ \ \frac{\sin(\sqrt{1-m\tau}|\xi|t)}{\sqrt{1-m\tau}|\xi|}\mathrm{e}^{-m\tau^2|\xi|^2t}.
 \end{align}
 Comparing \eqref{W-01} with \eqref{W-02} straightforwardly, the propagation speed has been changed from $1$ into $\sqrt{1-m\tau}$, and the new diffusion part $\mathrm{e}^{-m\tau^2|\xi|^2t}$ appears due to $m>0$. This new discovery is the root for the phenomenon of the optimal decay estimates found in Remark \ref{Rem-Influence-01}.
\end{remark}
Our second task is to describe more detailed information of the acoustic velocity potential $\psi$, namely, the optimal error estimates (by subtracting the optimal leading term) for large time under the new condition $Q_{\psi_0,\psi_1,\psi_2}\neq0$, which is 
\begin{align}\label{Quantity}
	|Q_{\psi_0,\psi_1,\psi_2}|^2:=|M_{\psi_1+\tau\psi_2}|^2+\left|(\delta+2m\tau^2)P_{\psi_0-\tau^2\psi_2}+\frac{n\gamma_0}{2\sqrt{1-m\tau}}P_{\psi_1+\tau\psi_2}\right|^2+|P_{\psi_1+\tau\psi_2}|^2,
\end{align}
with initial data localizing in the function space 
\begin{align*}
\ml{B}_{\delta,\ell}:=\ml{A}_{\delta,\ell}\cap (L^1\times L^{1,1}\times L^{1,1}).
\end{align*}
We recall the weighted $L^1$ space via
\begin{align*}
	L^{1,1}:=\left\{f\in L^1 \ :\ \|f\|_{L^{1,1}}:=\int_{\mb{R}^n}(1+|x|)|f(x)|\mathrm{d}x<+\infty \right\}.
\end{align*}
It means that we require an additional $L^{1,1}$ regularity for the second and third data.
\begin{remark}
In the viscous memoryless case $m=0$ with $\delta>0$, our condition $Q_{\psi_0,\psi_1,\psi_2}\neq0$ improves the one in \cite[Proposition 2.8 and Theorem 2.1]{Chen-Takeda=2023} because we do not need to consider the possibility for $P_{\psi_0-\tau^2\psi_2}\neq0$ in \eqref{Quantity}. Instead, we should check the chance of the sum
\begin{align*}
P_{\psi_0-\tau^2\psi_2}+\frac{(4\tau-\delta)\delta n}{16}P_{\psi_1+\tau\psi_2}\neq0
\end{align*}
for the viscous MGT equation \eqref{MGT} with $\delta>0$.
\end{remark}

\begin{theorem}\label{Thm-Optimal-Leading}
Suppose that initial data $(\psi_0,\psi_1,\psi_2)\in\ml{B}_{\delta,\ell}$ with $\ell>\frac{n}{2}$ and $\delta\geqslant0$. Then, the solution to the MGT equations \eqref{MGT-Memory} with the memory \eqref{Memory-Kernel} satisfies the following optimal error estimates:
\begin{align}\label{Est-Leading}
t^{-\frac{n}{4}}|Q_{\psi_0,\psi_1,\psi_2}|\lesssim\|\psi(t,\cdot)-\psi^{(1,p)}(t,\cdot)\|_{L^2}\lesssim t^{-\frac{n}{4}} \|(\psi_0,\psi_1,\psi_2)\|_{\ml{B}_{\delta,\ell}}
\end{align}
for large time $t\gg1$, provided that $Q_{\psi_0,\psi_1,\psi_2}\neq0$. Furthermore, the further error estimate holds
\begin{align}\label{Est-Further-Prof}
\|\psi(t,\cdot)-\psi^{(1,p)}(t,\cdot)-\psi^{(2,p)}(t,\cdot)\|_{L^2}=o(t^{-\frac{n}{4}})
\end{align}
for large time $t\gg1$.
\end{theorem}
\begin{remark}
		By subtracting the optimal leading term $\psi^{(1,p)}$ in the $L^2$ framework, the optimal estimates in Theorem \ref{Thm-Optimal-Est} can be improved by $t^{\frac{n}{4}}\ml{D}_n(t)$ for large time. Particularly, this improvement is optimal provided that the new condition $Q_{\psi_0,\psi_1,\psi_2}\neq 0$ holds. Due to an  improvement on the decay rate by subtracting the additional function $\psi^{(2,p)}$ in \eqref{Est-Further-Prof}, we may explain it as a second  profile of the solution for large time.
\end{remark}

\subsection{Result and discussion on the inviscid limits}
$\ \ \ \ $We now turn to the inviscid limits for the Cauchy problem \eqref{MGT-Memory} with respect to small diffusivity of sound $0<\delta\ll 1$. 
 Let us introduce the $L^{\infty}$ energy term for a function $u=u(t,x)$ as follows:
\begin{align*}
\|u(t,\cdot)\|_{\mathrm{E}_{\infty}}:=\left\|\left(|D|u_t+\frac{\delta+2\tau}{2\tau(\delta+\tau)}|D|u\right)(t,\cdot)\right\|_{L^{\infty}}+\left\|\left(u_{tt}+\frac{\delta+2\tau}{2\tau(\delta+\tau)}u_t\right)(t,\cdot)\right\|_{L^{\infty}}.
\end{align*}
Notice that $\frac{1}{2\tau}<\frac{\delta+2\tau}{2\tau(\delta+\tau)}<\frac{1}{\tau}$ for any $\delta,\tau>0$. We recall the Bessel potential space $H^s_1=\langle D\rangle^{-s}L^1$ and the Riesz potential space of negative order $\dot{H}^{-s}_1=|D|^sL^1$ with $s>0$, where $|D|^s$ and $\langle D\rangle^{-s}$ stand for the pseudo-differential operators with the symbols $|\xi|^s$ and $\langle \xi\rangle^{-s}$, respectively.

\begin{theorem}\label{Thm-Inviscid-Limit}
Let $0<\delta\ll 1$.  We denote $\psi^{\delta>0}=\psi^{\delta>0}(t,x)$ and $\psi^{\delta=0}=\psi^{\delta=0}(t,x)$ to be the solutions for the viscous model \eqref{MGT-Memory} with $\delta>0$ and the inviscid model \eqref{MGT-Memory} with $\delta=0$, respectively, owing the matching initial data, namely,
\begin{align*}
	\partial_t^k\psi^{\delta>0}(0,x)=\partial_t^k\psi(0,x)=\partial_t^k\psi^{\delta=0}(0,x)\ \ \mbox{with}\ \ k=0,1,2.
\end{align*}
\begin{itemize}
	\item By assuming $(\psi_0,\psi_1,\psi_2)\in(H^{s_0+1}_1)^2\times H^{s_0}_1$ with $s_0>n+2$, the energy term satisfies the following global (in time) convergence:
	\begin{align*}
		\sup\limits_{t\in[0,+\infty)}\|\psi^{\delta>0}(t,\cdot)-\psi^{\delta=0}(t,\cdot)\|_{\mathrm{E}_{\infty}}\leqslant C\sqrt{\delta}\|(\psi_0,\psi_1,\psi_2)\|_{(H^{s_0+1}_1)^2\times H^{s_0}_1}.
	\end{align*}
\item By assuming $(\psi_0,\psi_1,\psi_2)\in(H^{s_1+1}_1\cap \dot{H}^{-1}_1)^2\times (H^{s_1}_1\cap \dot{H}^{-1}_1)$ with $s_1>n+1$, the acoustic velocity potential satisfies the following local (in time) convergence:
\begin{align*}
	\sup\limits_{t\in[0,T]}\|\psi^{\delta>0}(t,\cdot)-\psi^{\delta=0}(t,\cdot)\|_{L^{\infty}}\leqslant C_{T}\sqrt{\delta}\|(\psi_0,\psi_1,\psi_2)\|_{(H^{s_1+1}_1\cap \dot{H}^{-1}_1)^2\times (H^{s_1}_1\cap \dot{H}^{-1}_1)}.
\end{align*}
\end{itemize}
Hereinbefore, the constant $C$ is independent of $\delta$.
\end{theorem}
\begin{remark}
The last uniform estimates suggest that $\psi^{\delta>0}(t,\cdot)$ converges to $\psi^{\delta=0}(t,\cdot)$ in $\mathrm{E}_{\infty}$ globally and in $L^{\infty}$ locally with the rate $\sqrt{\delta}$ of convergence, which exactly coincide with the formal WKB expansion derived in Proposition \ref{Prop-Formal-WKB}.
\end{remark}
\begin{remark}
The restriction of local (in time) convergence for the acoustic velocity potential comes from the memory kernel. To be specific, we actually may get the  global (in time) memory-weighted convergence as follows:
\begin{align*}
\sup\limits_{t\in[0,+\infty)}\left(\sqrt{g^{\tau}(t)}\|\psi^{\delta>0}(t,\cdot)-\psi^{\delta=0}(t,\cdot)\|_{L^{\infty}}\right)\leqslant C\sqrt{\delta}\|(\psi_0,\psi_1,\psi_2)\|_{(H^{s_1+1}_1\cap \dot{H}^{-1}_1)^2\times (H^{s_1}_1\cap \dot{H}^{-1}_1)},
\end{align*}
where the constant $C$ is independent of $\delta$ and $t$. Therefore, the constant $C_{T}$ in Theorem \ref{Thm-Inviscid-Limit} can be explicitly written by $C_T=\mathrm{e}^{\frac{T}{2\tau}}$. The question of global (in time) convergence for the acoustic velocity potential in the $L^{\infty}$ framework seems challenging due to the result in Theorem \ref{Thm-Optimal-Est} that the solutions $\psi^{\delta>0}$ and $\psi^{\delta=0}$ fulfill different kinds of decay estimates (no-loss type versus regularity-loss type).
\end{remark}

\section{Large time asymptotic behavior of the solution}\label{Section-MGT-Large-Time}\setcounter{equation}{0}
$\ \ \ \ $This section devotes to some large time qualitative properties, including the optimal growth/decay estimates and the optimal leading term in the $L^2$ framework, of the solution to the Cauchy problem \eqref{MGT-Memory} by using the WKB analysis and the Fourier analysis.

 Before carrying out our deduction, motivated by \cite[Section 1]{Dell-Pata=2017}, the first insight is transferring the integro-differential equation \eqref{MGT-Memory} into the following fourth order (in time) evolution equation:
\begin{align}\label{Fourth-PDE}
\begin{cases}
\ml{L}_{\mathrm{MGT-M}}^{4-\mathrm{th}}(\partial_t,|D|)\psi=0,&x\in\mb{R}^n,\ t>0,\\
(\psi,\psi_t,\psi_{tt},\psi_{ttt})(0,x)=(\psi_0,\psi_1,\psi_2,\psi_3)(x),&x\in\mb{R}^n,
\end{cases}
\end{align}
with the differential operator
\begin{align*}
\ml{L}_{\mathrm{MGT-M}}^{(4-\mathrm{th})}(\partial_t,|D|):=\tau^2\partial_t^4+2\tau\partial_t^3+\big(I+\tau(\delta+\tau)|D|^2\big)\partial_t^2+(\delta+2\tau)|D|^2\partial_t+(1-m\tau)|D|^2,
\end{align*}
where we acted the operator $I+\tau\partial_t$ equipped the identity operator $I$ to the equation in the Cauchy problem \eqref{MGT-Memory} because of
\begin{align}\label{reduce-g}
(I+\tau\partial_t)(g^{\tau}\ast\Delta \psi)=m\tau\Delta\psi
\end{align}
thanks to the formula \eqref{Memory-Kernel}.
Here, the last initial data is given by
\begin{align}\label{Four-Data}
\psi_3(x)=\frac{1}{\tau}\big(\Delta\psi_0(x)+(\delta+\tau)\Delta\psi_1(x)-\psi_2(x)\big).
\end{align}
Benefit from the determined initial condition \eqref{Four-Data}, the solution to the Cauchy problem \eqref{MGT-Memory} is uniquely given by
\begin{align*}
\psi(t,x)=K_0(t,|D|)\psi_0(x)+K_1(t,|D|)\psi_1(x)+K_2(t,|D|)\psi_2(x).
\end{align*}
We later will analyze asymptotic behavior of these kernels $K_j(t,|D|)$ with $j=0,1,2$ according to the fundamental solutions to the fourth order PDEs \eqref{Fourth-PDE}.

By applying the partial Fourier transform in regard to spatial variables to the Cauchy problem \eqref{Fourth-PDE}, one immediately obtains
\begin{align}\label{Fourth-Fourier}
	\begin{cases}
		\ml{L}_{\mathrm{MGT-M}}^{(4-\mathrm{th})}(\partial_t,|\xi|)\widehat{\psi}=0,&\xi\in\mb{R}^n,\ t>0,\\
		(\widehat{\psi},\widehat{\psi}_t,\widehat{\psi}_{tt},\widehat{\psi}_{ttt})(0,\xi)=(\widehat{\psi}_0,\widehat{\psi}_1,\widehat{\psi}_2,\widehat{\psi}_3)(\xi),&\xi\in\mb{R}^n,
	\end{cases}
\end{align}
whose corresponding characteristic equation is 
\begin{align}\label{Quartic-Eq}
	\tau^2\lambda^4+2\tau\lambda^3+\big(1+\tau(\delta+\tau)|\xi|^2\big)\lambda^2+(\delta+2\tau)|\xi|^2\lambda+(1-m\tau)|\xi|^2=0.
\end{align}

As preparations, we take the next zones in the Fourier space:
\begin{align*}
	\ml{Z}_{\intt}(\varepsilon_0):=\{ |\xi|\leqslant\varepsilon_0\ll1\},\ \
	\ml{Z}_{\bdd}(\varepsilon_0,N_0):=\{\varepsilon_0\leqslant |\xi|\leqslant N_0\},\ \  
	\ml{Z}_{\extt}(N_0):=\{|\xi|\geqslant N_0\gg1\}.
\end{align*}
The cut-off functions $\chi_{\intt}(\xi),\chi_{\bdd}(\xi),\chi_{\extt}(\xi)\in \mathcal{C}^{\infty}$ having supports in their corresponding zones $\ml{Z}_{\intt}(\varepsilon_0)$, $\ml{Z}_{\bdd}(\varepsilon_0/2,2N_0)$ and $\ml{Z}_{\extt}(N_0)$, respectively, fulfilling $\chi_{\bdd}(\xi)=1-\chi_{\intt}(\xi)-\chi_{\extt}(\xi)$. We next will study asymptotic behaviors in these zones.

\subsection{Asymptotic analysis for small frequencies}\label{Sub-Sec-Asym-Small}
$\ \ \ \ $Before stating some pointwise estimates of the solution $\widehat{\psi}$ to the fourth order differential equation \eqref{Fourth-Fourier} for small frequencies $\xi\in\ml{Z}_{\intt}(\varepsilon_0)$, let us prepare the partial Fourier transform for $J_k$ with $k=0,1,2$ as follows:
\begin{align*}
	\widehat{J}_0:=\frac{\sin(\sqrt{1-m\tau}|\xi|t)}{\sqrt{1-m\tau}|\xi|}\mathrm{e}^{-\frac{1}{2}(\delta+2m\tau^2)|\xi|^2t}\ \ \mbox{and}\ \ \widehat{J}_1:=\cos(\sqrt{1-m\tau}|\xi|t)\mathrm{e}^{-\frac{1}{2}(\delta+2m\tau^2)|\xi|^2t},
\end{align*}
as well as
\begin{align*}
	\widehat{J}_2:=\frac{\cos(\sqrt{1-m\tau}|\xi|t)}{\sqrt{1-m\tau}}\gamma_0|\xi|^2t\mathrm{e}^{-\frac{1}{2}(\delta+2m\tau^2)|\xi|^2t}=\frac{\gamma_0}{\sqrt{1-m\tau}}|\xi|^2t\widehat{J}_1,
\end{align*}
with the constant $\gamma_0$ defined in \eqref{gamma-0}. Remark that the last functions are well-defined because of the condition $m\tau<1$ in \eqref{Memory-Kernel}. Indeed, the dissipative effect from the previous functions is mainly generated by two factors of the exponential function
\begin{align*}
\mathrm{e}^{-\frac{1}{2}(\delta+2m\tau^2)|\xi|^2t}=\mathrm{e}^{-\frac{1}{2}\delta|\xi|^2t}\times\mathrm{e}^{-m\tau^2|\xi|^2t},
\end{align*}
in which the first exponential factor $\mathrm{e}^{-\frac{1}{2}\delta|\xi|^2t}$ comes from the viscous term ($\delta>0$) of the model, and the second exponential factor $\mathrm{e}^{-m\tau^2|\xi|^2t}$ comes from the memory term ($m>0$) of the model. In other words, even for inviscid fluids with $\delta=0$, the exponential factor in the sense of the Gaussian kernel still works due to the memory term from hereditary fluids.
\begin{prop}\label{Prop-Small-Freq}
	The solution $\widehat{\psi}$ to the Cauchy problem \eqref{Fourth-Fourier} satisfies the following sharp pointwise estimates in the Fourier space:
	\begin{align*}
		\chi_{\intt}(\xi)|\widehat{\psi}|&\lesssim\chi_{\intt}(\xi)\left(1+\frac{|\sin(\sqrt{1-m\tau}|\xi|t)|}{|\xi|}\right)\mathrm{e}^{-c|\xi|^2t}\left(|\widehat{\psi}_0|+|\widehat{\psi}_1|+|\widehat{\psi}_2|\right),
	\end{align*}
	and 
	\begin{align*}	
		\chi_{\intt}(\xi)\left|\widehat{\psi}-\widehat{J}_0(\widehat{\psi}_1+\tau\widehat{\psi}_2)\right|&\lesssim\chi_{\intt}(\xi)\mathrm{e}^{-c|\xi|^2t}\left(|\widehat{\psi}_0|+|\widehat{\psi}_1|+|\widehat{\psi}_2|\right),\\
		\chi_{\intt}(\xi)\left|\widehat{\psi}-(\widehat{J}_0+\widehat{J}_2)(\widehat{\psi}_1+\tau\widehat{\psi}_2)-\widehat{J}_1(\widehat{\psi}_0-\tau^2\widehat{\psi}_2)\right|&\lesssim\chi_{\intt}(\xi)|\xi|\mathrm{e}^{-c|\xi|^2t}\left(|\widehat{\psi}_0|+|\widehat{\psi}_1|+|\widehat{\psi}_2|\right).
	\end{align*}
\end{prop}
\begin{proof}
First of all, due to the condition $1-m\tau>0$, the discriminant $\triangle_{\mathrm{dis}}$ of the quartic equation \eqref{Quartic-Eq} satisfies
\begin{align*}
\triangle_{\mathrm{dis}}=-16m\tau^5(1-m\tau)|\xi|^4+O(|\xi|^6)<0 \ \ \mbox{when}\ \ \xi\in\ml{Z}_{\intt}(\varepsilon_0),
\end{align*}
which leads to two distinct real roots $\lambda_{1,2}$ and two complex conjugate non-real roots $\lambda_{3,4}$. Precisely, thanks to Taylor's expansions with respect to $|\xi|$, the characteristic roots have the following asymptotic expansions with suitable orders:
\begin{align*}
\lambda_{1,2}=-\frac{1}{\tau}\pm\sqrt{m\tau}|\xi|+\frac{1}{2}(2m\tau^2-\delta)|\xi|^2+O(|\xi|^3),\ \ \lambda_{3,4}=\lambda_{\mathrm{R}}\pm i\lambda_{\mathrm{I}}
\end{align*}
carrying the real and imaginary parts, respectively,
\begin{align*}
\lambda_{\mathrm{R}}=-\frac{1}{2}(\delta+2m\tau^2)|\xi|^2+O(|\xi|^4)\ \ \mbox{and}\ \ \lambda_{\mathrm{I}}=\sqrt{1-m\tau}|\xi|+\gamma_0|\xi|^3+O(|\xi|^4),
\end{align*}
 where the constant $\gamma_0$ is defined in \eqref{gamma-0}. Note that the suitable order asymptotic expansions of these roots contribute to the correct constructions of $\widehat{J}_k$ with $k=0,1,2$, particularly,
 \begin{align*}
 \lambda_2-\lambda_1=-2\sqrt{m\tau}|\xi|+O(|\xi|^3)\ \ \mbox{and}\ \ \lambda_{\mathrm{I}}-\sqrt{1-m\tau}|\xi|=\gamma_0|\xi|^3+O(|\xi|^4).
 \end{align*}

Therefore, benefit from the conjugate roots, the unique solution to the Cauchy problem \eqref{Fourth-Fourier} is given by
\begin{align*}
	\chi_{\intt}(\xi)\widehat{\psi}=\chi_{\intt}(\xi)\left\{d_1\mathrm{e}^{\lambda_1t}+d_2\mathrm{e}^{\lambda_2t}+\mathrm{e}^{\lambda_{\mathrm{R}}t}[(d_3+d_4)\cos(\lambda_{\mathrm{I}}t)+i(d_3-d_4)\sin(\lambda_{\mathrm{I}}t)]\right\},
\end{align*}
where the coefficients $d_j$ of the exponential functions $\mathrm{e}^{\lambda_jt}$ are determined according to
\begin{align}\label{linear-argebra}
	\underbrace{\begin{pmatrix}
			1 & 1 &1 &1\\
			\lambda_1 & \lambda_2 & \lambda_3 & \lambda_4\\
			\lambda_1^2 & \lambda_2^2 & \lambda_3^2 & \lambda_4^2\\
			\lambda_1^3 & \lambda_2^3 & \lambda_3^3 & \lambda_4^3\\
	\end{pmatrix}}_{=:\mb{V}}
	\begin{pmatrix}
		d_1 \\
		d_2 \\
		d_3 \\
		d_4\\
	\end{pmatrix}
	=\underbrace{\begin{pmatrix}
			\widehat{\psi}_0 \\
			\widehat{\psi}_1 \\
			\widehat{\psi}_2 \\
			-\frac{1}{\tau}|\xi|^2\widehat{\psi}_0-\frac{\delta+\tau}{\tau}|\xi|^2\widehat{\psi}_1-\frac{1}{\tau}\widehat{\psi}_2\\
	\end{pmatrix}}_{=:\mb{D}}.
\end{align}
The determinant of this Vandermonde matrix $\mb{V}$ is expanded as follows:
\begin{align*}
\det(\mb{V})=\prod\limits_{1\leqslant j<k\leqslant 4}(\lambda_k-\lambda_j)&=-2i\lambda_{\mathrm{I}}\big(\lambda_4^2-(\lambda_1+\lambda_2)\lambda_4+\lambda_1\lambda_2\big)\big(\lambda_3^2-(\lambda_1+\lambda_2)\lambda_3+\lambda_1\lambda_2\big)\\
&\quad\times\big(-2\sqrt{m\tau}|\xi|+O(|\xi|^3)\big)\\
&=\frac{4i\sqrt{m\tau(1-m\tau)}}{\tau^4}|\xi|^2+O(|\xi|^4)\ \ \mbox{when}\ \ \xi\in\ml{Z}_{\intt}(\varepsilon_0).
\end{align*}
With the aid of the  well-known Cramer rule in the system of linear equations \eqref{linear-argebra}, the coefficients $d_j$ are expressed via
\begin{align}\label{dj}
	d_j=\frac{\det(\mb{V}_j)}{\det(\mb{V})}\ \ \mbox{for}\ \ j=1,\dots,4,
\end{align}
in which $\mb{V}_j$ denotes the matrix formed by replacing the corresponding $j$-th column of $\mb{V}$ by the column vector $\mb{D}$ defined in \eqref{linear-argebra}. Then, carrying out lengthy but straightforward calculations, the determinants of $\mb{V}_j$ can be shown explicitly, which are quite important to derive the refined representation of the solution, as follows:
\begin{align*}
	\det(\mb{V}_1)&=-2i\lambda_2\lambda_{\mathrm{I}}(\lambda_{\mathrm{R}}^2+\lambda_{\mathrm{I}}^2)\left((\lambda_{\mathrm{R}}-\lambda_2)^2+\lambda_{\mathrm{I}}^2\right)\widehat{\psi}_0\\
	&\quad+\left(2i\lambda_{\mathrm{I}}(\lambda_{\mathrm{R}}^2+\lambda_{\mathrm{I}}^2)^2+\lambda_2^2\lambda_4^2(\lambda_4-\lambda_2)-\lambda_2^2\lambda_3^2(\lambda_3-\lambda_2)\right)\widehat{\psi}_1\\
	&\quad+\left(-4i\lambda_{\mathrm{R}}\lambda_{\mathrm{I}}(\lambda_{\mathrm{R}}^2+\lambda_{\mathrm{I}}^2)-\lambda_2\lambda_4(\lambda_4^2-\lambda_2^2)+\lambda_2\lambda_3(\lambda_3^2-\lambda_2^2)\right)\widehat{\psi}_2\\
	&\quad+2i\lambda_{\mathrm{I}}\left((\lambda_{\mathrm{R}}-\lambda_2)^2+\lambda_{\mathrm{I}}^2\right)\left(-\tfrac{1}{\tau}|\xi|^2\widehat{\psi}_0-\tfrac{\delta+\tau}{\tau}|\xi|^2\widehat{\psi}_1-\tfrac{1}{\tau}\widehat{\psi}_2\right),
\end{align*}
\begin{align*}
	\det(\mb{V}_2)&=2i\lambda_1\lambda_{\mathrm{I}}(\lambda_{\mathrm{R}}^2+\lambda_{\mathrm{I}}^2)\left((\lambda_{\mathrm{R}}-\lambda_1)^2+\lambda_{\mathrm{I}}^2\right)\widehat{\psi}_0\\
	&\quad-\left(2i\lambda_{\mathrm{I}}(\lambda_{\mathrm{R}}^2+\lambda_{\mathrm{I}}^2)^2+\lambda_1^2\lambda_4^2(\lambda_4-\lambda_1)-\lambda_1^2\lambda_3^2(\lambda_3-\lambda_1)\right)\widehat{\psi}_1\\
	&\quad-\left(-4i\lambda_{\mathrm{R}}\lambda_{\mathrm{I}}(\lambda_{\mathrm{R}}^2+\lambda_{\mathrm{I}}^2)-\lambda_1\lambda_4(\lambda_4^2-\lambda_1^2)+\lambda_1\lambda_3(\lambda_3^2-\lambda_1^2)\right)\widehat{\psi}_2\\
	&\quad-2i\lambda_{\mathrm{I}}\left((\lambda_{\mathrm{R}}-\lambda_1)^2+\lambda_{\mathrm{I}}^2\right)\left(-\tfrac{1}{\tau}|\xi|^2\widehat{\psi}_0-\tfrac{\delta+\tau}{\tau}|\xi|^2\widehat{\psi}_1-\tfrac{1}{\tau}\widehat{\psi}_2\right),
\end{align*}
\begin{align*}
	\det(\mb{V}_3)&=\lambda_1\lambda_2\lambda_4(\lambda_4-\lambda_2)(\lambda_4-\lambda_1)(\lambda_2-\lambda_1)\widehat{\psi}_0\\
	&\quad-\left(\lambda_2^2\lambda_4^2(\lambda_4-\lambda_2)-\lambda_1^2\lambda_4^2(\lambda_4-\lambda_1)+\lambda_1^2\lambda_2^2(\lambda_2-\lambda_1)\right)\widehat{\psi}_1\\
	&\quad+\left(\lambda_2\lambda_4(\lambda_4^2-\lambda_2^2)-\lambda_1\lambda_4(\lambda_4^2-\lambda_1^2)+\lambda_1\lambda_2(\lambda_2^2-\lambda_1^2)\right)\widehat{\psi}_2\\
	&\quad-(\lambda_4-\lambda_2)(\lambda_4-\lambda_1)(\lambda_2-\lambda_1)\left(-\tfrac{1}{\tau}|\xi|^2\widehat{\psi}_0-\tfrac{\delta+\tau}{\tau}|\xi|^2\widehat{\psi}_1-\tfrac{1}{\tau}\widehat{\psi}_2\right),
\end{align*}
\begin{align*}
	\det(\mb{V}_4)&=-\lambda_1\lambda_2\lambda_3(\lambda_3-\lambda_2)(\lambda_3-\lambda_1)(\lambda_2-\lambda_1)\widehat{\psi}_0\\
	&\quad+\left(\lambda_2^2\lambda_3^2(\lambda_3-\lambda_2)-\lambda_1^2\lambda_3^2(\lambda_3-\lambda_1)+\lambda_1^2\lambda_2^2(\lambda_2-\lambda_1)\right)\widehat{\psi}_1\\
	&\quad-\left(\lambda_2\lambda_3(\lambda_3^2-\lambda_2^2)-\lambda_1\lambda_3(\lambda_3^2-\lambda_1^2)+\lambda_1\lambda_2(\lambda_2^2-\lambda_1^2)\right)\widehat{\psi}_2\\
	&\quad+(\lambda_3-\lambda_2)(\lambda_3-\lambda_1)(\lambda_2-\lambda_1)\left(-\tfrac{1}{\tau}|\xi|^2\widehat{\psi}_0-\tfrac{\delta+\tau}{\tau}|\xi|^2\widehat{\psi}_1-\tfrac{1}{\tau}\widehat{\psi}_2\right).
\end{align*}
Summing up the last determinants of $\mb{V}_j$ and the coefficients $d_j$ defined in \eqref{dj}, we claim the refined representation of the solution 
\begin{align}
	\chi_{\intt}(\xi)\widehat{\psi}&=\chi_{\intt}(\xi)\left\{\frac{\det(\mb{V}_1)}{\det(\mb{V})}\mathrm{e}^{\lambda_1t}+\frac{\det(\mb{V}_2)}{\det(\mb{V})}\mathrm{e}^{\lambda_2t}+\cos(\lambda_{\mathrm{I}}t)\frac{\det(\mb{V}_3)+\det(\mb{V}_4)}{\det(\mb{V})}\mathrm{e}^{\lambda_{\mathrm{R}}t}\right.\notag\\
	&\quad\quad\quad\quad\quad\left.+i\sin(\lambda_{\mathrm{I}}t)\frac{\det(\mb{V}_3)-\det(\mb{V}_4)}{\det(\mb{V})}\mathrm{e}^{\lambda_{\mathrm{R}}t}\right\}.\label{Rep-wide-psi}
\end{align}
Note that the last formula may describe the interplay between oscillations from $\sin(\lambda_{\mathrm{I}}t)$ and non-summable singularities from $\chi_{\intt}(\xi)[\det(\mb{V})]^{-1}$, which allows us to derive the optimal $L^2$ estimates in lower dimensions.

According to the asymptotic expansions of characteristic roots at the beginning of this section and the key cancellations in the coefficients of initial data $\widehat{\psi}_2$ localizing in $\det(\mb{V}_1)$ as well as $\det(\mb{V}_2)$, namely,
\begin{align*}
\lambda_{k}^2(\lambda_4-\lambda_3)-i\frac{2}{\tau}\lambda_{\mathrm{I}}\lambda_{k}^2=\mp i\frac{2}{\tau^2}\sqrt{m\tau(1-m\tau)}|\xi|^2+O(|\xi|^3)\ \ \mbox{with}\ \ k=1,2,
\end{align*}
we derive the following estimate:
\begin{align*}
\chi_{\intt}(\xi)\big(|\det(\mb{V}_1)|+|\det(\mb{V}_2)|\big)\lesssim\chi_{\intt}(\xi)\left(|\xi|^3|\widehat{\psi}_0|+|\xi|^3|\widehat{\psi}_1|+|\xi|^2|\widehat{\psi}_2|\right),
\end{align*}
which implies an exponential decay upper bound
\begin{align}\label{Solution-Est-01}
\chi_{\intt}(\xi)\left|\frac{\det(\mb{V}_1)}{\det(\mb{V})}\mathrm{e}^{\lambda_1t}+\frac{\det(\mb{V}_2)}{\det(\mb{V})}\mathrm{e}^{\lambda_2t}\right|\lesssim\chi_{\intt}(\xi)\mathrm{e}^{-ct}\left(|\widehat{\psi}_0|+|\widehat{\psi}_1|+|\widehat{\psi}_2|\right).
\end{align}
For this reason, the leading term and the higher order profiles of $\chi_{\intt}(\xi)\widehat{\psi}$ localize in the remaining parts of \eqref{Rep-wide-psi}. Concerning the cosine part, a direct calculation shows
\begin{align*}
\chi_{\intt}(\xi)|\det(\mb{V}_3)+\det(\mb{V}_4)|\lesssim\chi_{\intt}(\xi)\left(|\xi|^2|\widehat{\psi}_0|+|\xi|^4|\widehat{\psi}_1|+|\xi|^2|\widehat{\psi}_2|\right),
\end{align*}
whose dominant factor, i.e. the lowest order term with respect to the frequency $|\xi|$, is taken by
\begin{align*}
\det(\mb{V}_+)&:=-2i(\lambda_2-\lambda_1)\lambda_{\mathrm{I}}\lambda_1^2\lambda_2^2\widehat{\psi}_0+2i(\lambda_2-\lambda_1)\lambda_{\mathrm{I}}(\lambda_1^2-\lambda_1\lambda_2+\lambda_2^2)\widehat{\psi}_2\\
&\ =\left(\frac{4i\sqrt{m\tau(1-m\tau)}}{\tau^4}|\xi|^2+O(|\xi|^4)\right)\widehat{\psi}_0-\left(\frac{4i\sqrt{m\tau(1-m\tau)}}{\tau^2}|\xi|^2+O(|\xi|^4)\right)\widehat{\psi}_2
\end{align*}
when $\xi\in\ml{Z}_{\intt}(\varepsilon_0)$.
In other words, the refined estimate with an additional factor $|\xi|^2$ holds
\begin{align}\label{Error-Est-01}
\chi_{\intt}(\xi)|\det(\mb{V}_3)+\det(\mb{V}_4)-\det(\mb{V}_+)|\lesssim\chi_{\intt}(\xi)|\xi|^3\left(|\widehat{\psi}_0|+|\widehat{\psi}_1|+|\widehat{\psi}_2|\right).
\end{align}
It shows that
\begin{align}\label{Solution-Est-02}
\chi_{\intt}(\xi)\left|\cos(\lambda_{\mathrm{I}}t)\frac{\det(\mb{V}_3)+\det(\mb{V}_4)}{\det(\mb{V})}\mathrm{e}^{\lambda_{\mathrm{R}}t}\right|\lesssim\chi_{\intt}(\xi)\mathrm{e}^{-c|\xi|^2t}\left(|\widehat{\psi}_0|+|\widehat{\psi}_1|+|\widehat{\psi}_2|\right).
\end{align}
Finally, concerning the sine part, by the same method as the above, we find
\begin{align*}
\chi_{\intt}(\xi)|\det(\mb{V}_3)-\det(\mb{V}_4)|\lesssim\chi_{\intt}(\xi)\left(|\xi|^3|\widehat{\psi}_0|+|\xi|\,|\widehat{\psi}_1|+|\xi|\,|\widehat{\psi}_2|\right),
\end{align*}
whose dominant factor is extracted as follows:
\begin{align*}
\det(\mb{V}_-)&:=-2\lambda_1^2\lambda_2^2(\lambda_2-\lambda_1)\widehat{\psi}_1+2\lambda_1\lambda_2(\lambda_2-\lambda_1)\left(\frac{1}{\tau}+\lambda_1+\lambda_2\right)\widehat{\psi}_2\\
&\ =\left(\frac{4\sqrt{m\tau}}{\tau^4}|\xi|+O(|\xi|^3)\right)\widehat{\psi}_1+\left(\frac{4\sqrt{m\tau}}{\tau^3}|\xi|+O(|\xi|^3)\right)\widehat{\psi}_2
\end{align*}
when $\xi\in\ml{Z}_{\intt}(\varepsilon_0)$. Then, the refined estimate with an additional factor $|\xi|^2$ holds
\begin{align}\label{Error-Est-02}
\chi_{\intt}(\xi)|\det(\mb{V}_3)-\det(\mb{V}_4)-\det(\mb{V}_-)|\lesssim\chi_{\intt}(\xi)|\xi|^3\left(|\widehat{\psi}_0|+|\widehat{\psi}_1|+|\widehat{\psi}_2|\right).
\end{align}
It implies that
\begin{align}\label{Solution-Est-03}
	&\chi_{\intt}(\xi)\left|i\sin(\lambda_{\mathrm{I}}t)\frac{\det(\mb{V}_3)-\det(\mb{V}_4)}{\det(\mb{V})}\mathrm{e}^{\lambda_{\mathrm{R}}t}\right|\notag\\
	&\lesssim\chi_{\intt}(\xi)\mathrm{e}^{-c|\xi|^2t}\left(1+\frac{|\sin(\sqrt{1-m\tau}|\xi|t)|}{|\xi|}\right)\left(|\widehat{\psi}_0|+|\widehat{\psi}_1|+|\widehat{\psi}_2|\right),
\end{align}
where we used
\begin{align}
\chi_{\intt}(\xi)|\sin(\lambda_{\mathrm{I}}t)-\sin(\sqrt{1-m\tau}|\xi|t)|&=\chi_{\intt}(\xi)\left|\gamma_0|\xi|^3t+O(|\xi|^4)t\right||\cos(\eta_0t)|\notag\\
&\lesssim\chi_{\intt}(\xi)|\xi|^3t\label{Error-01}
\end{align}
from the mean value theorem with a $|\xi|$-dependent function $\eta_0=\eta_0(|\xi|)$ valued between $\sqrt{1-m\tau}|\xi|$ and $\lambda_{\mathrm{I}}$, moreover,
\begin{align*}
\chi_{\intt}(\xi)|\sin(\lambda_{\mathrm{I}}t)|\mathrm{e}^{-c|\xi|^2t}&\lesssim\chi_{\intt}(\xi)|\xi|^3t\mathrm{e}^{-c|\xi|^2t}+\chi_{\intt}(\xi)|\sin(\sqrt{1-m\tau}|\xi|t)|\mathrm{e}^{-c|\xi|^2t}\\
&\lesssim\chi_{\intt}(\xi)\left(|\xi|+|\sin(\sqrt{1-m\tau}|\xi|t)|\right)\mathrm{e}^{-c|\xi|^2t}.
\end{align*}
Combining with the derived estimates \eqref{Solution-Est-01}, \eqref{Solution-Est-02} and \eqref{Solution-Est-03}, we arrive at our first estimate. 

With the aid of \eqref{Error-Est-02} and \eqref{Error-Est-01}, respectively, we may obtain
\begin{align}\label{Err-02}
\chi_{\intt}(\xi)\left|\widehat{\psi}-i\sin(\lambda_{\mathrm{I}}t)\frac{\det(\mb{V}_-)}{\det(\mb{V})}\mathrm{e}^{\lambda_{\mathrm{R}}t}\right|\lesssim\chi_{\intt}(\xi)\mathrm{e}^{-c|\xi|^2t}\left(|\widehat{\psi}_0|+|\widehat{\psi}_1|+|\widehat{\psi}_2|\right)
\end{align}
and
\begin{align}
&\chi_{\intt}(\xi)\left|\widehat{\psi}-i\sin(\lambda_{\mathrm{I}}t)\frac{\det(\mb{V}_-)}{\det(\mb{V})}\mathrm{e}^{\lambda_{\mathrm{R}}t}-\cos(\lambda_{\mathrm{I}}t)\frac{\det(\mb{V}_+)}{\det(\mb{V})}\mathrm{e}^{\lambda_{\mathrm{R}}t}\right|\notag\\
&\lesssim\chi_{\intt}(\xi)|\xi|\mathrm{e}^{-c|\xi|^2t}\left(|\widehat{\psi}_0|+|\widehat{\psi}_1|+|\widehat{\psi}_2|\right).\label{Err-03}
\end{align}
To achieve our desire refined estimates, let us focus on suitable approximations for the next terms:
\begin{align*}
\widehat{G}_0:=i\sin(\lambda_{\mathrm{I}}t)\frac{\det(\mb{V}_-)}{\det(\mb{V})}\mathrm{e}^{\lambda_{\mathrm{R}}t}\ \ \mbox{and}\ \ \widehat{G}_1:=\cos(\lambda_{\mathrm{I}}t)\frac{\det(\mb{V}_+)}{\det(\mb{V})}\mathrm{e}^{\lambda_{\mathrm{R}}t}
\end{align*}
for small frequencies. According to the facts that
\begin{align*}
\chi_{\intt}(\xi)\left|\frac{i\det(\mb{V}_-)\sqrt{1-m\tau}|\xi|-\det(\mb{V})(\widehat{\psi}_1+\tau\widehat{\psi}_2)}{\det(\mb{V})\sqrt{1-m\tau}|\xi|}\right|\lesssim\chi_{\intt}(\xi)|\xi|,
\end{align*}
and
\begin{align*}
\chi_{\intt}(\xi)\left|\mathrm{e}^{\lambda_{\mathrm{R}}t}-\mathrm{e}^{-\frac{1}{2}(\delta+2m\tau^2)|\xi|^2t}\right|&\lesssim\chi_{\intt}(\xi)|\xi|^4t\mathrm{e}^{-\frac{1}{2}(\delta+2m\tau^2)|\xi|^2t}\left|\int_0^1\mathrm{e}^{O(|\xi|^4)t\eta}\mathrm{d}\eta\right|\\
&\lesssim\chi_{\intt}(\xi)|\xi|^2\mathrm{e}^{-c|\xi|^2t},
\end{align*}
as well as the error estimate \eqref{Error-01}, we may get
\begin{align}
\chi_{\intt}(\xi)\left|\widehat{G}_0-\widehat{J}_0(\widehat{\psi}_1+\tau\widehat{\psi}_2)\right|&\lesssim\chi_{\intt}(\xi)|\sin(\lambda_{\mathrm{I}}t)|\left|\frac{i\det(\mb{V}_-)}{\det(\mb{V})}-\frac{\widehat{\psi}_1+\tau\widehat{\psi}_2}{\sqrt{1-m\tau}|\xi|}\right|\mathrm{e}^{\lambda_{\mathrm{R}}t}\notag\\
&\quad+\chi_{\intt}(\xi)\frac{|\sin(\lambda_{\mathrm{I}}t)|}{\sqrt{1-m\tau}|\xi|}\left|\mathrm{e}^{\lambda_{\mathrm{R}}t}-\mathrm{e}^{-\frac{1}{2}(\delta+2m\tau^2)|\xi|^2t}\right|\left(|\widehat{\psi}_1|+|\widehat{\psi}_2|\right)\notag\\
&\quad+\chi_{\intt}(\xi)\frac{|\sin(\lambda_{\mathrm{I}}t)-\sin(\sqrt{1-m\tau}|\xi|t)|}{\sqrt{1-m\tau}|\xi|}\mathrm{e}^{-\frac{1}{2}(\delta+2m\tau^2)|\xi|^2t}\left(|\widehat{\psi}_1|+|\widehat{\psi}_2|\right)\notag\\
&\lesssim\chi_{\intt}(\xi)\mathrm{e}^{-c|\xi|^2t}\left(|\widehat{\psi}_1|+|\widehat{\psi}_2|\right),\label{Err-04}
\end{align} 
whose worst part is generated from the error estimate for the sine functions. To explore a further refined estimate, especially for the worst error part, Taylor's expansion implies
\begin{align*}
\sin(\lambda_{\mathrm{I}}t)-\sin(\sqrt{1-m\tau}|\xi|t)-\cos(\sqrt{1-m\tau}|\xi|t)\gamma_0|\xi|^3t=O(|\xi|^6)t^2.
\end{align*}
It leads to the next estimate with an additional factor $|\xi|$:
\begin{align}
\chi_{\intt}(\xi)\left|\widehat{G}_0-(\widehat{J}_0+\widehat{J}_2)(\widehat{\psi}_1+\tau\widehat{\psi}_2)\right|\lesssim\chi_{\intt}(\xi)|\xi|\mathrm{e}^{-c|\xi|^2t}\left(|\widehat{\psi}_1|+|\widehat{\psi}_2|\right).\label{Err-05}
\end{align}
Analogously, it holds
\begin{align}
\chi_{\intt}(\xi)\left|\widehat{G}_1-\widehat{J}_1(\widehat{\psi}_0-\tau^2\widehat{\psi}_2)\right|\lesssim\chi_{\intt}(\xi)|\xi|\mathrm{e}^{-c|\xi|^2t}\left(|\widehat{\psi}_0|+|\widehat{\psi}_2|\right),\label{Err-06}
\end{align}
where we employed
\begin{align*}
\chi_{\intt}(\xi)|\cos(\lambda_{\mathrm{I}}t)-\cos(\sqrt{1-m\tau}|\xi|t)|&\lesssim\chi_{\intt}(\xi)|\xi|^3t,\\
\chi_{\intt}(\xi)\left|\frac{\det(\mb{V}_+)-\det(\mb{V})(\widehat{\psi}_0-\tau^2\widehat{\psi}_2)}{\det(\mb{V})}\right|&\lesssim\chi_{\intt}(\xi)|\xi|^2.
\end{align*}
Finally, summarizing \eqref{Err-02} with \eqref{Err-04}, \eqref{Err-03} with \eqref{Err-05} and \eqref{Err-06}, via the triangle inequality, we immediately conclude the our result.
\end{proof}
\begin{remark}
For the viscous memoryless case $m=0$ with $\delta>0$, i.e. the viscous MGT equation \eqref{MGT}, our estimates derived in Proposition \ref{Prop-Small-Freq} exactly coincide with those in \cite[Proposition 2.1]{Chen-Ikehata=2021} and \cite[Proposition 2.4]{Chen-Takeda=2023}, where we also correct the misprint of $\widehat{\ml{N}}_0(t,|\xi|)$ with an additional minus, i.e. the correct function $\widehat{J}_2$, in \cite[Proposition 2.4]{Chen-Takeda=2023}.
\end{remark}

\subsection{Asymptotic analysis for large frequencies}\label{Sub-Sec-Asym-Large}
$\ \ \ \ $Let us now turn to large frequencies $\xi\in\ml{Z}_{\extt}(N_0)$, particularly, some pointwise estimates of the solution $\widehat{\psi}$ to the fourth order differential equation \eqref{Fourth-Fourier}.
\begin{prop}\label{Prop-Large-Freq}
	The solution $\widehat{\psi}$ to the Cauchy problem \eqref{Fourth-Fourier} satisfies the following pointwise estimates in the Fourier space:
	\begin{align*}
	\chi_{\extt}(\xi)|\widehat{\psi}|\lesssim \chi_{\extt}(\xi)\left(|\widehat{\psi}_0|+|\xi|^{-1}|\widehat{\psi}_1|+|\xi|^{-2}|\widehat{\psi}_2|\right)\times\begin{cases}
	\mathrm{e}^{-ct}&\mbox{when}\ \ \delta>0,\\
	\mathrm{e}^{-c|\xi|^{-2}t}&\mbox{when}\ \ \delta=0.
	\end{cases}
	\end{align*}
\end{prop}
\begin{proof}
Thanks to $m>0$, the discriminant of the quartic equation \eqref{Quartic-Eq} fulfills
\begin{align*}
\triangle_{\mathrm{dis}}=-4\tau^5(\delta+\tau)^3\big(\delta^2+4m\tau^2(\delta+\tau)\big)|\xi|^{10}+O(|\xi|^{8})<0\ \ \mbox{when}\ \ \xi\in\ml{Z}_{\extt}(N_0).
\end{align*}
In other words, there exist two distinct real roots $\lambda_{1,2}$ and two complex conjugate non-real roots $\lambda_{3,4}$ as follows:
\begin{align*}
\lambda_{1,2}=\gamma_1^{\pm}+O(|\xi|^{-2}),\ \ \lambda_{3,4}=\lambda_{\mathrm{R}}\pm i\lambda_{\mathrm{I}}
\end{align*}
with a pair of negative constants
\begin{align*}
\gamma_1^{\pm}:=\frac{-(\delta+2\tau)\pm\sqrt{\delta^2+4m\tau^2(\delta+\tau)}}{2\tau(\delta+\tau)},
\end{align*}
carrying the real and imaginary parts, respectively,
\begin{align}\label{Real-Part}
\begin{cases}
\displaystyle{\lambda_{\mathrm{R}}=-\frac{\delta}{2\tau(\delta+\tau)}+O(|\xi|^{-2})\ \ \mbox{and}\ \ \lambda_{\mathrm{I}}=\sqrt{\frac{\delta+\tau}{\tau}}|\xi|+O(|\xi|^{-1})}&\mbox{when}\ \ \delta>0,\\[1em]
\displaystyle{\lambda_{\mathrm{R}}=-\frac{m}{\tau^2}|\xi|^{-2}+O(|\xi|^{-4})\ \ \mbox{and}\ \ \lambda_{\mathrm{I}}=|\xi|+\frac{m}{2\tau}|\xi|^{-1}+O(|\xi|^{-3})}&\mbox{when}\ \ \delta=0.
\end{cases}
\end{align}

\begin{remark}
As we have explained in the introduction, the condition $m\tau<1$ in \eqref{Memory-Kernel} always holds physically due to the typically small relaxation parameters $m$ and $\tau$. If one considers the counterfact $m\tau>1$ mathematically, we may find $\gamma_1^+>0$ since $\delta^2+4m\tau^2(\delta+\tau)>(\delta+2\tau)^2$. As a consequence, we immediately conclude an exponential instability for the MGT equation in hereditary fluids \eqref{MGT-Memory} when $m\tau>1$.
\end{remark}

In order to estimate each term in the same representation \eqref{Rep-wide-psi} with the localized function $\chi_{\extt}(\xi)$, we may follow the same ways as those in \eqref{linear-argebra}, where we will omit their details for the sake of similarity and briefness.  We prepare
\begin{align*}
\det(\mb{V})=\frac{2i(\delta+\tau)^2}{\tau^3}\sqrt{\frac{\delta^2}{\tau(\delta+\tau)}+4m\tau}\,|\xi|^5+O(|\xi|^4)\ \ \mbox{when}\ \ \xi\in\ml{Z}_{\extt}(N_0).
\end{align*}
Concerning any $\delta\geqslant0$, by performing lengthy but straightforward calculations similarly to those in Subsection \ref{Sub-Sec-Asym-Small}, we may arrive at
 	\begin{align*}
 	\chi_{\extt}(\xi)\big(|\det(\mb{V}_1)|+|\det(\mb{V}_2)|\big)\lesssim\chi_{\extt}(\xi)\left(|\xi|^5|\widehat{\psi}_0|+|\xi|^3|\widehat{\psi}_1|+|\xi|^3|\widehat{\psi}_2|\right)
 \end{align*}
 due to the cancellation in the coefficient of $\widehat{\psi}_1$ that
 \begin{align*}
 	2i\lambda_{\mathrm{I}}^5-\frac{\delta+\tau}{\tau}2i\lambda_{\mathrm{I}}^3|\xi|^2=2i\lambda_{\mathrm{I}}^3\left(\lambda_{\mathrm{I}}^2-\frac{\delta+\tau}{\tau}|\xi|^2\right)=O(|\xi|^3),
 \end{align*}
and
\begin{align*}
	\chi_{\extt}(\xi)|\det(\mb{V}_3)+\det(\mb{V}_4)|&\lesssim\chi_{\extt}(\xi)\left(|\xi|^3|\widehat{\psi}_0|+|\xi|^3|\widehat{\psi}_1|+|\xi|^3|\widehat{\psi}_2|\right),\\
	\chi_{\extt}(\xi)|\det(\mb{V}_3)-\det(\mb{V}_4)|&\lesssim\chi_{\extt}(\xi)\left(|\xi|^4|\widehat{\psi}_0|+|\xi|^4|\widehat{\psi}_1|+|\xi|^2|\widehat{\psi}_2|\right).
\end{align*}
Although the real part $\lambda_{\mathrm{R}}$ has different behavior between the viscous case $\delta>0$ and the inviscid case $\delta=0$, this term does not influence on the dominant parts of $\det(\mb{V}_j)$ with $j=1,\dots,4$. The combination of the last estimates in the representation \eqref{Rep-wide-psi} with the localized function $\chi_{\extt}(\xi)$ completes the proof of our proposition.
\end{proof}
\begin{remark}
The threshold condition for decay properties as $\xi\in\ml{Z}_{\extt}(N_0)$ is addressed by $\delta=0$. To be specific, an exponential decay property holds when $\delta>0$; a polynomial decay property in the 2-regularity-loss type holds when $\delta=0$; an exponential instability holds when $\delta<0$, according to the real parts expanded in \eqref{Real-Part}.
\end{remark}

\subsection{Exponential stability for bounded frequencies}
$\ \ \ \ $Finally, we investigate the situation  $\xi\in\ml{Z}_{\bdd}(\varepsilon_0,N_0)$. To demonstrate negativities for the real parts of characteristic roots, by fixing an index $j\in\{1,\dots,4\}$, let us assume a contrario that $\lambda_j=i\tilde{\lambda}_j$ with  $\tilde{\lambda}_j\in\mb{R}\backslash\{0\}$. From the quartic equation \eqref{Quartic-Eq}, it fulfills the next system:
\begin{align*}
	\begin{cases}
		\tau^2\tilde{\lambda}_j^4-\big(1+\tau(\delta+\tau)|\xi|^2\big)\tilde{\lambda}_j^2+(1-m\tau)|\xi|^2=0,\\
		-2\tau \tilde{\lambda}_j^3+(\delta+2\tau)|\xi|^2\tilde{\lambda}_j=0.
	\end{cases}
\end{align*}
Because of $\tilde{\lambda}_j\neq0$, the second equation shows $\tilde{\lambda}_j^2=\frac{\delta+2\tau}{2\tau}|\xi|^2$. Plugging this value into the first equation implies
\begin{align*}
	0=-\frac{\delta(\delta+2\tau)}{4}|\xi|^4-\frac{2m\tau^2+\delta}{2\tau}|\xi|^2<0,
\end{align*}
which yields a contradiction. Moreover, let us recall $\mathrm{Re}\,\lambda_k<0$ with $k=1,\dots,4$ for $\xi\in\ml{Z}_{\intt}(\varepsilon_0)\cup\ml{Z}_{\extt}(N_0)$ in the last two subsections. So, the continuity of characteristic roots with respect to $|\xi|$ implies that the real parts of all characteristic roots are negative. The next exponential decay estimate holds.
\begin{prop}\label{Prop-Bdd-Freq}
	The solution $\widehat{\psi}$ to the Cauchy problem \eqref{Fourth-Fourier} satisfies the following pointwise estimate in the Fourier space:
\begin{align*}
	\chi_{\bdd}(\xi)|\widehat{\psi}|\lesssim \chi_{\bdd}(\xi)\mathrm{e}^{-ct}\left(|\widehat{\psi}_0|+|\widehat{\psi}_1|+|\widehat{\psi}_2|\right).
\end{align*}
\end{prop}

\subsection{Optimal growth/decay estimates: Proof of Theorem \ref{Thm-Optimal-Est}}
$\ \ \ \ $Before deriving the optimal $L^2$ estimates of the solution $\psi$ to the main problem \eqref{MGT-Memory}, let us recall the next useful lemma obtained by \cite{Ikehata=2014,Ikehata-Ono=2017}.
\begin{lemma}\label{Lem-Ike}
Let $c_0,c_1>0$. The Fourier multipliers fulfill the following optimal estimates:
\begin{align*}
\left\|\cos(c_0|\xi|t)\mathrm{e}^{-c_1|\xi|^2t}\right\|_{L^2}\simeq t^{-\frac{n}{4}}\ \ \mbox{and}\ \ \left\|\frac{\sin(c_0|\xi|t)}{|\xi|}\mathrm{e}^{-c_1|\xi|^2t}\right\|_{L^2}\simeq\ml{D}_n(t)
\end{align*}
for large time $t\gg1$, where the time-dependent function $\ml{D}_n(t)$ is defined in \eqref{Dnt}.
\end{lemma}

Firstly, by applying the last lemma and Proposition \ref{Prop-Small-Freq}, we have
\begin{align*}
\|\chi_{\intt}(D)\psi(t,\cdot)\|_{L^2}&\lesssim\left(\left\|\chi_{\intt}(\xi)\mathrm{e}^{-c|\xi|^2t}\right\|_{L^2}+\left\|\chi_{\intt}(\xi)\frac{\sin(c_0|\xi|t)}{|\xi|}\mathrm{e}^{-c|\xi|^2t}\right\|_{L^2}\right)\left\|\left(\widehat{\psi}_0,\widehat{\psi}_1,\widehat{\psi}_2\right)\right\|_{(L^{\infty})^3}\\
&\lesssim\ml{D}_n(1+t)\|(\psi_0,\psi_1,\psi_2)\|_{(L^1)^3}
\end{align*}
with $c_0=\sqrt{1-m\tau}>0$, where the Hausdorff-Young inequality and the Plancherel theorem were used. Due to the sharp estimate
\begin{align*}
\left\|\chi_{\extt}(\xi)\mathrm{e}^{-c|\xi|^{-2}t}\widehat{g}(\xi)\right\|_{L^2}&\lesssim\left\|\chi_{\extt}(\xi)|\xi|^{-\ell}\mathrm{e}^{-c|\xi|^{-2}t}\right\|_{L^{\infty}}\|\langle\xi\rangle^{\ell}\widehat{g}(\xi)\|_{L^2}\\
&\lesssim(1+t)^{-\frac{\ell}{2}}\|g\|_{H^{\ell}}
\end{align*}
for a given function $g=g(x)\in H^{\ell}$ with $\ell\geqslant0$, from Proposition \ref{Prop-Large-Freq}, one derives
\begin{align*}
\|\chi_{\extt}(D)\psi(t,\cdot)\|_{L^2}\lesssim\begin{cases}
\mathrm{e}^{-ct}\|(\psi_0,\psi_1,\psi_2)\|_{(L^2)^3}&\mbox{when}\ \ \delta>0,\\
(1+t)^{-\frac{\ell}{2}}\|(\psi_0,\psi_1,\psi_2)\|_{H^{\ell}\times H^{\ell-1}\times H^{\ell-2}}&\mbox{when}\ \ \delta=0.
\end{cases}
\end{align*}
Additionally from the exponential decay estimate stated in Proposition \ref{Prop-Bdd-Freq}, the following upper bound estimate holds:
\begin{align}\label{General-Est}
\|\psi(t,\cdot)\|_{L^2}&\lesssim\|\chi_{\intt}(D)\psi(t,\cdot)\|_{L^2}+\|\chi_{\bdd}(D)\psi(t,\cdot)\|_{L^2}+\|\chi_{\extt}(D)\psi(t,\cdot)\|_{L^2}\notag\\
&\lesssim\begin{cases}
	\ml{D}_n(1+t)\|(\psi_0,\psi_1,\psi_2)\|_{\ml{A}_{\delta,\ell}}&\mbox{when}\ \ \delta>0,\\
	\max\left\{\ml{D}_n(1+t),(1+t)^{-\frac{\ell}{2}}\right\}\|(\psi_0,\psi_1,\psi_2)\|_{\ml{A}_{0,\ell}}&\mbox{when}\ \ \delta=0,
\end{cases}
\end{align}
with any $\ell\geqslant0$. Letting $\ell\geqslant\max\{\frac{n}{2}-1,0\}$ to fix the dominant part in the above maximum when $\delta=0$, we get the desire upper bound estimate in \eqref{Est-Itself}.

To guarantee the sharpness of the above estimate, we now turn to its lower bound for large time. By the same way as preceding parts of the text associated with Proposition \ref{Prop-Small-Freq}, one finds
\begin{align*}
&\left\|\psi(t,\cdot)-J_0(t,|D|)\big(\psi_1(\cdot)+\tau\psi_2(\cdot)\big)\right\|_{L^2}\\
&\lesssim\left\|\chi_{\intt}(\xi)\mathrm{e}^{-c|\xi|^2t}\right\|_{L^2}\|(\psi_0,\psi_1,\psi_2)\|_{(L^1)^3}+\left\|\big(1-\chi_{\intt}(\xi)\big)\widehat{\psi}(t,\xi)\right\|_{L^2}+\mathrm{e}^{-ct}\|(\psi_1,\psi_2)\|_{(L^2)^2}\\
&\lesssim\begin{cases}
	(1+t)^{-\frac{n}{4}}\|(\psi_0,\psi_1,\psi_2)\|_{\ml{A}_{\delta,\ell}}&\mbox{when}\ \ \delta>0,\\
	\max\left\{(1+t)^{-\frac{n}{4}},(1+t)^{-\frac{\ell}{2}}\right\}\|(\psi_0,\psi_1,\psi_2)\|_{\ml{A}_{0,\ell}}&\mbox{when}\ \ \delta=0,
\end{cases}
\end{align*}
with any $\ell\geqslant0$. Then, with the aid of the mean value theorem
\begin{align*}
|J_0(t,x-y)-J_0(t,x)|\lesssim|y|\,|\nabla J_0(t,x-\eta_1y)|
\end{align*}
with $\eta_1\in(0,1)$, we may separate the integral into two parts such that
\begin{align}\label{Est-01}
&\left\|J_0(t,|D|)\big(\psi_1(\cdot)+\tau\psi_2(\cdot)\big)-J_0(t,\cdot)P_{\psi_1+\tau\psi_2}\right\|_{L^2}\notag\\
&\lesssim\left\|\int_{|y|\leqslant t^{\frac{1}{16}}}\big(J_0(t,\cdot-y)-J_0(t,\cdot)\big)\big(\psi_1(y)+\tau\psi_2(y)\big)\mathrm{d}y\right\|_{L^2}\notag\\
&\quad+\left\|\int_{|y|\geqslant t^{\frac{1}{16}}}\big(|J_0(t,\cdot-y)|+|J_0(t,\cdot)|\big)|\psi_1(y)+\tau\psi_2(y)|\mathrm{d}y\right\|_{L^2}\notag\\
&\lesssim t^{\frac{1}{16}}\|\,|\xi|\widehat{J}_0(t,|\xi|)\|_{L^2}\|(\psi_1,\psi_2)\|_{(L^1)^2}+\|\widehat{J}_0(t,|\xi|)\|_{L^2}\|\psi_1+\tau\psi_2\|_{L^1(|x|\geqslant t^{\frac{1}{16}})}\notag\\
&\lesssim t^{\frac{1}{16}-\frac{n}{4}}\|(\psi_1,\psi_2)\|_{(L^1)^2}+o\big(\ml{D}_n(t)\big)
\end{align}
for large time $t\gg1$, thanks to our assumption on the $L^1$ integrabilities of $\psi_1$ as well as $\psi_2$, namely,
\begin{align*}
\lim\limits_{t\to+\infty}\int_{|y|\geqslant t^{\frac{1}{16}}}|\psi_1(y)+\tau\psi_2(y)|\mathrm{d}y=0.
\end{align*}
 Hence, the triangle inequality shows
\begin{align*}
\|\psi(t,\cdot)-J_0(t,\cdot)P_{\psi_1+\tau\psi_2}\|_{L^2}=o\big(\ml{D}_n(t)\big)
\end{align*}
by choosing $\ell>\max\{\frac{n}{2}-1,0\}$ leading to $t^{-\frac{\ell}{2}}<\ml{D}_n(t)$ for large time $t\gg1$. Eventually, let us employ the Minkowski inequality and combine the last estimate with Lemma \ref{Lem-Ike} to arrive at
\begin{align*}
\|\psi(t,\cdot)\|_{L^2}&\gtrsim\|J_0(t,\cdot)\|_{L^2}|P_{\psi_1+\tau\psi_2}|-\|\psi(t,\cdot)-J_0(t,\cdot)P_{\psi_1+\tau\psi_2}\|_{L^2}\\
&\gtrsim\ml{D}_n(t)|P_{\psi_1+\tau\psi_2}|-o\big(\ml{D}_n(t)\big)
\end{align*}
for large time $t\gg1$. It gives immediately our desire result.

\subsection{Optimal leading term: Proof of Theorem \ref{Thm-Optimal-Leading}}
$\ \ \ \ $According to the hypothesis $\psi_1+\tau\psi_2\in L^{1,1}$ and \cite[Lemma 2.2]{Ikehata=2014}, we may deduce
\begin{align*}
|\widehat{\psi}_1+\tau\widehat{\psi}_2-P_{\psi_1+\tau\psi_2}|\lesssim|\xi|\,\|(\psi_1,\psi_2)\|_{(L^{1,1})^2},
\end{align*}
so that the estimate \eqref{Est-01} can be improved by
\begin{align*}
\left\|J_0(t,|D|)\big(\psi_1(\cdot)+\tau\psi_2(\cdot)\big)-J_0(t,\cdot)P_{\psi_1+\tau\psi_2}\right\|_{L^2}&\lesssim\|\,|\xi|\widehat{J}_0(t,|\xi|)\|_{L^2}\|(\psi_1,\psi_2)\|_{(L^{1,1})^2}\\
&\lesssim (1+t)^{-\frac{n}{4}}\|(\psi_1,\psi_2)\|_{(L^{1,1})^2}.
\end{align*}
Hence, for any $\ell\geqslant0$, we have
\begin{align*}
\|\psi(t,\cdot)-J_0(t,\cdot)P_{\psi_1+\tau\psi_2}\|_{L^2}\lesssim\begin{cases}
	(1+t)^{-\frac{n}{4}}\|(\psi_0,\psi_1,\psi_2)\|_{\ml{B}_{\delta,\ell}}&\mbox{when}\ \ \delta>0,\\
	\max\left\{(1+t)^{-\frac{n}{4}},(1+t)^{-\frac{\ell}{2}}\right\}\|(\psi_0,\psi_1,\psi_2)\|_{\ml{B}_{0,\ell}}&\mbox{when}\ \ \delta=0,
\end{cases}
\end{align*}
where the relation $\ml{B}_{\delta,\ell}\subset\ml{A}_{\delta,\ell}$ was used. Choosing $\ell\geqslant\frac{n}{2}$ to fix the dominant part in the above maximum when $\delta=0$, one may claim the desire upper bound estimate in \eqref{Est-Leading}.

Let us focus on the further error estimate by the next suitable decomposition:
\begin{align*}
\psi-\psi^{(1,p)}-\psi^{(2,p)}=\ml{E}_1+\ml{E}_2+\ml{E}_3+\ml{E}_4,
\end{align*}
whose components $\ml{E}_k=\ml{E}_k(t,x)$ with $k=1,\dots,4$ are chosen by
\begin{align*}
\ml{E}_1&:=\psi-\big(J_0(t,|D|)+J_2(t,|D|)\big)(\psi_1+\tau\psi_2)-J_1(t,|D|)(\psi_0-\tau^2\psi_2),\\
\ml{E}_2&:=J_0(t,|D|)(\psi_1+\tau\psi_2)-J_0(t,x)P_{\psi_1+\tau\psi_2}+\nabla J_0(t,x)\circ M_{\psi_1+\tau\psi_2},\\
\ml{E}_3&:=J_2(t,|D|)(\psi_1+\tau\psi_2)-J_2(t,x)P_{\psi_1+\tau\psi_2},\\
\ml{E}_4&:=J_1(t,|D|)(\psi_0-\tau^2\psi_2)-J_1(t,x)P_{\psi_0-\tau^2\psi_2}.
\end{align*}
Recalling the third pointwise estimate in Proposition \ref{Prop-Small-Freq}, we know
\begin{align*}
\|\chi_{\intt}(\xi)\widehat{\ml{E}}_1(t,\xi)\|_{L^2}&\lesssim\left\|\chi_{\intt}(\xi)|\xi|\mathrm{e}^{-c|\xi|^2t}\right\|_{L^2}\|(\psi_0,\psi_1,\psi_2)\|_{(L^1)^3}\\
&\lesssim (1+t)^{-\frac{1}{2}-\frac{n}{4}}\|(\psi_0,\psi_1,\psi_2)\|_{(L^1)^3}.
\end{align*}
Concerning bounded and large frequencies, it holds
\begin{align*}
\left\|\big(\chi_{\bdd}(\xi)+\chi_{\extt}(\xi)\big)\widehat{\ml{E}}_1(t,\xi)\right\|_{L^2}\lesssim\begin{cases}
	\mathrm{e}^{-ct}\|(\psi_0,\psi_1,\psi_2)\|_{(L^2)^3}&\mbox{when}\ \ \delta>0,\\
	(1+t)^{-\frac{\ell}{2}}\|(\psi_0,\psi_1,\psi_2)\|_{H^{\ell}\times H^{\ell-1}\times H^{\ell-2}}&\mbox{when}\ \ \delta=0,
\end{cases}
\end{align*}
so that
\begin{align*}
\|\ml{E}_1(t,\cdot)\|_{L^2}\lesssim\begin{cases}
	(1+t)^{-\frac{1}{2}-\frac{n}{4}}\|(\psi_0,\psi_1,\psi_2)\|_{\ml{A}_{\delta,\ell}}&\mbox{when}\ \ \delta>0,\\
	\max\left\{(1+t)^{-\frac{1}{2}-\frac{n}{4}},(1+t)^{-\frac{\ell}{2}}\right\}\|(\psi_0,\psi_1,\psi_2)\|_{\ml{A}_{0,\ell}}&\mbox{when}\ \ \delta=0,
\end{cases} 
\end{align*}
with any $\ell\geqslant0$. For the third and fourth error terms, analogously to \eqref{Est-01}, we are able to deduce
\begin{align*}
\|\ml{E}_3(t,\cdot)\|_{L^2}+\|\ml{E}_4(t,\cdot)\|_{L^2}=o(t^{-\frac{n}{4}})
\end{align*}
for large time $t\gg1$, because of the assumptions $\psi_1+\tau\psi_2\in L^1$ and $\psi_0-\tau^2\psi_2\in L^1$. In order to treat the final error term $\ml{E}_2$, we make use of Taylor's expansion of second order
\begin{align*}
|J_0(t,x-y)-J_0(t,x)+y\circ \nabla J_0(t,x)|\lesssim|y|^2|\nabla^2J_0(t,x-\eta_2y)|
\end{align*}
with $\eta_2\in(0,1)$. Moreover, the separation $\ml{E}_2=\ml{E}_{2,1}+\ml{E}_{2,2}$ holds with $\ml{E}_{2,k}=\ml{E}_{2,k}(t,x)$ such that
\begin{align*}
\ml{E}_{2,1}(t,x)&:=\int_{|y|\leqslant t^{\frac{1}{16}}}\big(J_0(t,x-y)-J_0(t,x)+y\circ\nabla J_0(t,x)\big)\big(\psi_1(y)+\tau\psi_2(y)\big)\mathrm{d}y,\\
\ml{E}_{2,2}(t,x)&:=\int_{|y|\geqslant t^{\frac{1}{16}}}\big(J_0(t,x-y)-J_0(t,x)\big)\big(\psi_1(y)+\tau\psi_2(y)\big)\mathrm{d}y\\
&\quad+\int_{|y|\geqslant t^{\frac{1}{16}}}y\circ \nabla J_0(t,x)\big(\psi_1(y)+\tau\psi_2(y)\big)\mathrm{d}y.
\end{align*}
Therefore, one arrives at
\begin{align*}
\|\ml{E}_{2,1}(t,\cdot)\|_{L^2}&\lesssim t^{\frac{1}{8}}\|\,|\xi|^2\widehat{J}_0(t,|\xi|)\|_{L^2}\|(\psi_1,\psi_2)\|_{(L^1)^2}\lesssim t^{-\frac{3}{8}-\frac{n}{4}}\|(\psi_1,\psi_2)\|_{(L^1)^2},\\
\|\ml{E}_{2,2}(t,\cdot)\|_{L^2}&\lesssim\|\,|\xi|\widehat{J}_0(t,|\xi|)\|_{L^2}\|(|x|\psi_1,|x|\psi_2)\|_{(L^1(|x|\geqslant t^{\frac{1}{16}}))^2}=o(t^{-\frac{n}{4}}),
\end{align*}
for large time $t\gg1$, thanks to the additional weighted $L^1$ assumption $\psi_1,\psi_2\in L^{1,1}$. From the relation $\ml{B}_{\delta,\ell}\subset\ml{A}_{\delta,\ell}$ again, the sum of last derived estimates gives
\begin{align*}
\sum\limits_{k=1,\dots,4}\|\ml{E}_k(t,\cdot)\|_{L^2}\lesssim \begin{cases}
	o(t^{-\frac{n}{4}})\|(\psi_0,\psi_1,\psi_2)\|_{\ml{B}_{\delta,\ell}}&\mbox{when}\ \ \delta>0,\\
	\max\left\{o(t^{-\frac{n}{4}}),(1+t)^{-\frac{\ell}{2}}\right\}\|(\psi_0,\psi_1,\psi_2)\|_{\ml{B}_{0,\ell}}&\mbox{when}\ \ \delta=0,
\end{cases}
\end{align*}
with any $\ell\geqslant0$. So, taking $\ell>\frac{n}{2}$, the further error estimate \eqref{Est-Further-Prof} is obtained.

In order to derive the optimal lower bound estimate for the error term
\begin{align*}
\|\psi(t,\cdot)-\psi^{(1,p)}(t,\cdot)\|_{L^2}\ \ \mbox{for large time}\ \ t\gg1,
\end{align*}
namely, to determine the optimal leading term, the crucial part is to estimate the second profile $\psi^{(2,p)}(t,\cdot)$ in the $L^2$ norm from the below side. To begin with our deduction, the partial Fourier transform of it is addressed by
\begin{align*}
\widehat{\psi}^{(2,p)}=-i\xi\widehat{J}_0\circ M_{\psi_1+\tau\psi_2}+\widehat{J}_1\left(P_{\psi_0-\tau^2\psi_2}+\frac{\gamma_0t}{\sqrt{1-m\tau}}|\xi|^2P_{\psi_1+\tau\psi_2}\right).
\end{align*}
Thus, our aim is separated into
\begin{align*}
\|\widehat{\psi}^{(2,p)}(t,\xi)\|_{L^2}^2=\int_{\mb{R}^n}\left(|\mathrm{Im}\,\widehat{\psi}^{(2,p)}(t,\xi)|^2+|\mathrm{Re}\,\widehat{\psi}^{(2,p)}(t,\xi)|^2\right)\mathrm{d}\xi=A_1(t)+A_2(t),
\end{align*}
in which we defined
\begin{align*}
A_1(t)&:=\int_{\mb{R}^n}|\widehat{J}_0(t,|\xi|)|^2|\xi\circ M_{\psi_1+\tau\psi_2}|^2\mathrm{d}\xi,\\
A_2(t)&:=\int_{\mb{R}^n}|\widehat{J}_1(t,|\xi|)|^2\left|P_{\psi_0-\tau^2\psi_2}+\frac{\gamma_0t}{\sqrt{1-m\tau}}|\xi|^2P_{\psi_1+\tau\psi_2}\right|^2\mathrm{d}\xi.
\end{align*}
For one thing, the use of polar coordinates shows
\begin{align}\label{A1}
A_1(t)=\int_0^{+\infty}r^{n+1}|\widehat{J}_0(t,r)|^2\mathrm{d}r\int_{|\omega|=1}|\omega\circ M_{\psi_1+\tau\psi_2}|^2\mathrm{d}\sigma_{\omega}\gtrsim t^{-\frac{n}{2}}|M_{\psi_1+\tau\psi_2}|^2
\end{align}
for large time $t\gg1$. For another, setting $\eta=rt^{\frac{1}{2}}$, we may separate $A_2(t)$ into three parts as follows:
\begin{align*}
A_2(t)&=C_n\int_0^{+\infty}r^{n-1}\left|P_{\psi_0-\tau^2\psi_2}+\frac{\gamma_0 t}{\sqrt{1-m\tau}}r^2P_{\psi_1+\tau\psi_2}\right|^2|\cos(\sqrt{1-m\tau} rt)|^2\,\mathrm{e}^{-(\delta+2m\tau^2)r^2t}\mathrm{d}r\\
&=C_n\int_0^{+\infty}r^{n-1}\left(|P_{\psi_0-\tau^2\psi_2}|^2+\frac{2\gamma_0t}{\sqrt{1-m\tau}}r^2P_{\psi_0-\tau^2\psi_2}P_{\psi_1+\tau\psi_2}+\frac{\gamma_0^2t^2}{1-m\tau}r^4|P_{\psi_1+\tau\psi_2}|^2 \right)\\
&\qquad\qquad\quad\times|\cos(\sqrt{1-m\tau}rt)|^2\,\mathrm{e}^{-(\delta+2m\tau^2)r^2t}\mathrm{d}r\\
&=C_nt^{-\frac{n}{2}}\left(|P_{\psi_0-\tau^2\psi_2}|^2A_3^{(0)}(t)+\frac{2\gamma_0}{\sqrt{1-m\tau}}P_{\psi_0-\tau^2\psi_2}P_{\psi_1+\tau\psi_2}A_3^{(1)}(t)+\frac{\gamma_0^2}{1-m\tau}|P_{\psi_1+\tau\psi_2}|^2A_3^{(2)}(t)\right)
\end{align*}
with a positive constant $C_n$, where we took the time-dependent function with $k=0,1,2$ such that
\begin{align*}
A_3^{(k)}(t):=\int_0^{+\infty}\eta^{n-1+2k}\left|\cos\big(\sqrt{(1-m\tau)t}\eta\big)\right|^2\mathrm{e}^{-(\delta+2m\tau^2)\eta^2}\mathrm{d}\eta.
\end{align*}
From the well-known fact $\cos^2y=\frac{1}{2}+\frac{1}{2}\cos(2y)$ and the Riemann-Lebesgue theorem, we derive
\begin{align*}
A_3^{(k)}(t)&=\frac{1}{2}\int_0^{+\infty}\eta^{n-1+2k}\mathrm{e}^{-(\delta+2m\tau^2)\eta^2}\mathrm{d}\eta+\frac{1}{2}\int_0^{+\infty}\eta^{n-1+2k}\cos\big(2\sqrt{(1-m\tau)t}\eta\big)\mathrm{e}^{-(\delta+2m\tau^2)\eta^2}\mathrm{d}\eta\\
&=\frac{1}{2}(\delta+2m\tau^2)^{-\frac{n}{2}-k}\int_0^{+\infty}\tilde{\eta}^{n-1+2k}\mathrm{e}^{-\tilde{\eta}^2}\mathrm{d}\tilde{\eta}+o(1)\\
&=\frac{1}{4}(\delta+2m\tau^2)^{-\frac{n}{2}-k}\,\Gamma\left(\frac{n}{2}+k\right)+o(1)
\end{align*}
for large time $t\gg1$,  where we considered the Gamma function
\begin{align*}
\Gamma(z):=\int_0^{+\infty}\eta^{z-1}\mathrm{e}^{-\eta}\mathrm{d}\eta=2\int_0^{+\infty}\tilde{\eta}^{2z-1}\mathrm{e}^{-\tilde{\eta}^2}\mathrm{d}\tilde{\eta}.
\end{align*}
It indicates
\begin{align*}
A_2(t)=\frac{C_n}{4}t^{-\frac{n}{2}}(\delta+2m\tau^2)^{-\frac{n}{2}-2}\widetilde{A}_2(t),
\end{align*}
where we use the property of the Gamma function, i.e. $\Gamma(z+1)=z\Gamma(z)$, and get
\begin{align*}
\widetilde{A}_2(t)&=(\delta+2m\tau^2)^2\,\Gamma\left(\frac{n}{2}\right)|P_{\psi_0-\tau^2\psi_2}|^2+\frac{2\gamma_0}{\sqrt{1-m\tau}}(\delta+2m\tau^2)\Gamma\left(\frac{n}{2}+1\right)P_{\psi_0-\tau^2\psi_2}P_{\psi_1+\tau\psi_2}\\
&\quad+\frac{\gamma_0^2}{1-m\tau}\Gamma\left(\frac{n}{2}+2\right)|P_{\psi_1+\tau\psi_2}|^2+o(1)\\
&=\Gamma\left(\frac{n}{2}\right)\left[\left|(\delta+2m\tau^2)P_{\psi_0-\tau^2\psi_2}+\frac{n\gamma_0}{2\sqrt{1-m\tau}}P_{\psi_1+\tau\psi_2}\right|^2+\frac{n\gamma_0^2}{2(1-m\tau)}|P_{\psi_1+\tau\psi_2}|^2\right]+o(1)
\end{align*}
for large time $t\gg1$, namely,
\begin{align}\label{A2}
A_2(t)\gtrsim t^{-\frac{n}{2}}\left[\left|(\delta+2m\tau^2)P_{\psi_0-\tau^2\psi_2}+\frac{n\gamma_0}{2\sqrt{1-m\tau}}P_{\psi_1+\tau\psi_2}\right|^2+|P_{\psi_1+\tau\psi_2}|^2 \right].
\end{align}
Recalling the quantity $Q_{\psi_0,\psi_1,\psi_2}$ defined according to \eqref{Quantity}, applying  the Minkowski inequality as well as \eqref{Est-Further-Prof} associated with \eqref{A1} and \eqref{A2}, we are able to claim
\begin{align*}
\|\psi(t,\cdot)-\psi^{(1,p)}(t,\cdot)\|_{L^2}&\gtrsim A_1(t)+A_2(t)-\|\psi(t,\cdot)-\psi^{(1,p)}(t,\cdot)-\psi^{(2,p)}(t,\cdot)\|_{L^2}\\
&\gtrsim t^{-\frac{n}{4}}|Q_{\psi_0,\psi_1,\psi_2}|-o(t^{-\frac{n}{4}})
\end{align*}
for large time $t\gg1$. It completes the lower bound estimate in \eqref{Est-Leading}.

\section{Inviscid limits of the solution in the $L^{\infty}$ framework}\setcounter{equation}{0}\label{Section-Inviscid}
\subsection{Formal expansion with the small diffusivity of sound}
$\ \ \ \ $Our main purpose of this section is to analyze some influence of the diffusivity of sound $0<\delta\ll 1$ on asymptotic behavior of the solution to the MGT equation \eqref{MGT-Memory} in hereditary fluids. That is to say, the viscous solution $\psi^{\delta>0}$ converges to the inviscid solution $\psi^{\delta=0}$ for the Cauchy problem \eqref{MGT-Memory} with the memory term \eqref{Memory-Kernel} in the $L^{\infty}$ framework.

 Strongly motivated by the Prandtl boundary layer theory and the multi-scale analysis, we look for an approximate viscous solution by a formal power series expansion in $\delta$ for the Cauchy problem \eqref{MGT-Memory} with $\delta>0$ via the WKB analysis, namely,
\begin{align}\label{Expansion-Solution-delta}
	\psi^{\delta>0}(t,x)=\sum\limits_{j\geqslant0}\delta^{\frac{j}{2}}\left(\psi^{I,j}(t,x)+\psi^{L,j}(t,z)\right)\ \ \mbox{with}\ \ z:=\frac{x}{\sqrt{\delta}}\in\mb{R}^n,
\end{align}
where all of the above terms in the representation are assumed to be smooth. Here,  $\psi^{I,j}=\psi^{I,j}(t,x)$ stands for the dominant profile for the $j$-th order expansion and $\psi^{L,j}=\psi^{L,j}(t,z)$ denotes the other profile decaying to zero as $|z|\to+\infty$. We will determine these profiles formally soon afterwards.

Then, we may plug the formal expansion \eqref{Expansion-Solution-delta} into the equation of the Cauchy problem \eqref{MGT-Memory} to arrive at
\begin{align*}
	0&=\tau\sum\limits_{j\geqslant0}\delta^{\frac{j}{2}}\left(\psi_{ttt}^{I,j}+\psi_{ttt}^{L,j}\right)+\sum\limits_{j\geqslant0}\delta^{\frac{j}{2}}\left(\psi_{tt}^{I,j}+\psi_{tt}^{L,j}\right)-\sum\limits_{j\geqslant0}\delta^{\frac{j}{2}}\left(\Delta\psi^{I,j}+\delta^{-1}\Delta_{z}\psi^{L,j}\right)\\
	&\quad-\sum\limits_{j\geqslant0}\delta^{\frac{j}{2}+1}\left(\Delta\psi_t^{I,j}+\delta^{-1}\Delta_z\psi_t^{L,j}\right)-\tau\sum\limits_{j\geqslant0}\delta^{\frac{j}{2}}\left(\Delta\psi_t^{I,j}+\delta^{-1}\Delta_z\psi_t^{L,j}\right)\\
	&\quad+g^{\tau}\ast\sum\limits_{j\geqslant0}\delta^{\frac{j}{2}}\left(\Delta\psi^{I,j}+\delta^{-1}\Delta_z\psi^{L,j}\right),
\end{align*}
in which the Laplacian in terms of $z$ is denoted by $\Delta_z:=\sum_{k=1}^n\partial_{z_k}^2$. Matching the terms with the size $\delta^{\frac{j}{2}}$, we claim
\begin{align*}
	0&=\tau\left(\psi_{ttt}^{I,j}+\psi_{ttt}^{L,j}\right)+\left(\psi_{tt}^{I,j}+\psi_{tt}^{L,j}\right)-\left(\Delta\psi^{I,j}+\Delta_z\psi^{L,j+2}\right)-\left(\Delta\psi_t^{I,j-2}+\Delta_z\psi^{L,j}_t\right)\\
	&\quad-\tau\left(\Delta\psi^{I,j}_t+\Delta_z\psi_t^{L,j+2}\right)+g^{\tau}\ast\left(\Delta\psi^{I,j}+\Delta_z\psi^{L,j+2}\right).
\end{align*}
Naturally, we always impose $\psi^{I,j}\equiv0\equiv\psi^{L,j}$ when $j<0$. Letting $\delta\downarrow0$, i.e. $|z|\to+\infty$, the profiles $\psi^{L,j}$, $\psi^{L,j+2}$ and their derivatives tend to zero so that
\begin{align}\label{Eq-psi-I}
	\tau\psi_{ttt}^{I,j}+\psi_{tt}^{I,j}-\Delta\psi^{I,j}-\tau\Delta\psi_t^{I,j}+g^{\tau}\ast\Delta\psi^{I,j}=\Delta\psi^{I,j-2}_t,
\end{align}
as well as
\begin{align}\label{Eq-psi-L}
	\Delta_z\left(\tau\psi_t^{L,j+2}+\psi^{L,j+2}-g^{\tau}\ast\psi^{L,j+2}\right)=\partial_t\left(\tau\psi_{tt}^{L,j}-\Delta_z\psi^{L,j}+\psi_t^{L,j}\right),
\end{align}
whose initial conditions will be determined later. The quantities to be known functions on the right-hand sides of the equations \eqref{Eq-psi-I} and \eqref{Eq-psi-L} can be fixed by the last orders matches.

Moreover, one notices
\begin{align*}
	\partial_t^k\psi^{\delta>0}(0,x)=\sum\limits_{j\geqslant0}\delta^{\frac{j}{2}}\left(\partial_t^k\psi^{I,j}(0,x)+\partial_t^k\psi^{L,j}(0,z)\right)\ \ \mbox{with}\ \ k=0,1,2.
\end{align*}
Due to the initial conditions of the Cauchy problem \eqref{MGT-Memory} being independent of $\delta$, we then may naturally consider $\partial_t^k\psi^{I,j}(0,x)\equiv0\equiv\partial_t^k\psi^{L,j}(0,z)$ when $j\geqslant 1$ with $k=0,1,2$. As a result, we are able to deduce
\begin{align}\label{Restriction-Initial-Condition}
	\partial_t^k\psi^{I,0}(0,x)+\partial_t^k\psi^{L,0}(0,z)=\psi^{\delta>0}_k(x)\ \ \mbox{with}\ \ k=0,1,2.
\end{align}

Let us take $j=0$ in \eqref{Eq-psi-I} to find
\begin{align*}
	\begin{cases}
		\tau\psi^{I,0}_{ttt}+\psi^{I,0}_{tt}-\Delta\psi^{I,0}-\tau\Delta\psi^{I,0}_t+g^{\tau}\ast\Delta\psi^{I,0}=0,&x\in\mb{R}^n,\ t>0,\\
		\psi^{I,0}(0,x)=\psi^{I,0}_0(x),\ \psi^{I,0}_t(0,x)=\psi^{I,0}_1(x),\ \psi^{I,0}_{tt}(0,x)=\psi^{I,0}_2(x),&x\in\mb{R}^n.
	\end{cases}
\end{align*}
Moreover, letting $j=-2$ in \eqref{Eq-psi-L} and recalling $\psi^{L,-2}\equiv0$, we next apply the operator $I+\tau\partial_t$ on the resultant. With the help of \eqref{reduce-g}, it results
\begin{align}\label{Layer-L-1}
	\begin{cases}
		\Delta_z\left(\tau^2\psi_{tt}^{L,0}+2\tau\psi^{L,0}_t+(1-m\tau)\psi^{L,0}\right)=0,&z\in\mb{R}^n,\ t>0,\\
		\psi^{L,0}(0,z)=\psi^{L,0}_0(z),\ \psi^{L,0}_t(0,z)=\psi^{L,0}_1(z),\ \psi^{L,0}_{tt}(0,z)=\psi^{L,0}_2(z),&z\in\mb{R}^n.
	\end{cases}
\end{align}
 Note that the initial conditions in the last two problems need to be determined by themselves carefully rather than artificial choices. For this reason, let us introduce the auxiliary function $\psi^{P,0}=\psi^{P,0}(t,x)$ such that
\begin{align}\label{Expression-P0}
	\psi^{P,0}(t,x):=\psi^{L,0}(t,z)+\psi^{I,0}_0(x)+t\psi^{I,0}_1(x)+\frac{1}{2}t^2\psi^{I,0}_2(x)\ \ \mbox{with}\ \ z=\frac{x}{\sqrt{\delta}}.
\end{align}
According to \eqref{Layer-L-1} and \eqref{Restriction-Initial-Condition}, it solves the following Cauchy problem:
\begin{align*}
	\begin{cases}
		\displaystyle{\Delta\left(\tau^2\psi^{P,0}_{tt}+2\tau\psi^{P,0}_t+(1-m\tau)\psi^{P,0}\right)=\Delta F,}&x\in\mb{R}^n,\ t>0,\\
		\psi^{P,0}(0,x)=\psi_0^{\delta>0}(x),\ \psi^{P,0}_t(0,x)=\psi^{\delta>0}_1(x),\ \psi^{P,0}_{tt}(0,x)=\psi^{\delta>0}_2(x),&x\in\mb{R}^n,
	\end{cases}
\end{align*}
whose source term $F=F(t,x)$ is given by
\begin{align*}
F(t,x)&=\left(\tau^2\psi_2^{I,0}(x)+2\tau\psi^{I,0}_1(x)+(1-m\tau)\psi^{I,0}_0(x)\right)\\
&\quad+t\left(2\tau\psi^{I,0}_2(x)+(1-m\tau)\psi^{I,0}_1(x)\right)+\frac{1-m\tau}{2}t^2\psi_2^{I,0}(x).
\end{align*}
A direct computation from the classical ODE theory shows that the solution to the last over-determined Cauchy problem with three initial conditions is uniquely solved by
\begin{align*}
\psi^{P,0}(t,x)=\psi^{I,0}_0(x)+t\psi^{I,0}_1(x)+\frac{1}{2}t^2\psi^{I,0}_2(x),
\end{align*}
in which we applied some comparisons for the coefficients of time-dependent functions. Hence, in the view of \eqref{Expression-P0}, we conclude $\psi^{L,0}\equiv0$ and $\psi^{I,0}\equiv\psi^{\delta=0}$ with initial data $\psi^{\delta=0}_k(x)=\psi^{\delta>0}_k(x)$ for $k=0,1,2$.

Summarizing the last statements, we claim the formal WKB expansion of the solution. 
\begin{prop}\label{Prop-Formal-WKB}
	The solution $\psi^{\delta>0}$ to the viscous MGT equation \eqref{MGT-Memory} with the small diffusivity of sound $0<\delta\ll 1$ formally has the next asymptotic expansions:
	\begin{align*}
		\psi^{\delta>0}(t,x)=\psi^{\delta=0}(t,x)+\sum\limits_{j\geqslant 0}\delta^{\frac{j+1}{2}}\left(\psi^{I,j+1}(t,x)+\psi^{L,j+1}(t,\tfrac{x}{\sqrt{\delta}})\right),
	\end{align*}
	where $\psi^{\delta=0}$ is the solution to the inviscid MGT equation \eqref{MGT-Memory} with the vanishing diffusivity of sound $\delta=0$, $\psi^{I,j+1}$ and $\psi^{L,j+1}$ solve the inhomogeneous integro-differential equations \eqref{Eq-psi-I} and \eqref{Eq-psi-L}, respectively, for any $j\geqslant0$.
\end{prop}

\subsection{Rigorous justification of the inviscid limits: Proof of Theorem \ref{Thm-Inviscid-Limit}}
$\ \ \ \ $According to Proposition \ref{Prop-Formal-WKB}, let us introduce the difference $u=u(t,x)$ between the viscous solution and the inviscid solution, i.e. $u=\psi^{\delta>0}-\psi^{\delta=0}$, fulfilling the next inhomogeneous MGT equation with vanishing initial data:
\begin{align}\label{Inhomo-MGT-Memory}
\begin{cases}
\tau u_{ttt}+u_{tt}-\Delta u-(\delta+\tau)\Delta u_t+g^{\tau}\ast\Delta u=-\delta\Delta \psi_t^{\delta=0},&x\in\mb{R}^n,\ t>0,\\
u(0,x)=0,\ u_t(0,x)=0, \ u_{tt}(0,x)=0,&x\in\mb{R}^n,
\end{cases}
\end{align}
where we employed the matching conditions $\psi^{\delta>0}_k(x)=\psi^{\delta=0}_k(x)$ with $k=0,1,2$.

We are motivated by \cite{Lasiecka-Wang-2016} to construct the first energy in the Fourier space
\begin{align*}
\mathbb{E}_1[\hat{u}]&:=\tau|\hat{u}_{tt}|^2+(\delta+\tau)|\xi|^2|\hat{u}_t|^2+2|\xi|^2\mathrm{Re}(\hat{u}\bar{\hat{u}}_t)+\frac{|\xi|^2}{\tau}\int_0^tg^{\tau}(t-\eta)|\hat{u}(t,\xi)-\hat{u}(\eta,\xi)|^2\mathrm{d}\eta\\
&\ \quad+g^{\tau}(t)|\xi|^2|\hat{u}|^2-2|\xi|^2\int_0^tg^{\tau}(t-\eta)\mathrm{Re}\big(\hat{u}(\eta,\xi)\bar{\hat{u}}_t(t,\xi)\big)\mathrm{d}\eta.
\end{align*}
Thanks to the properties $g^{\tau}(0)=m$ and $g^{\tau}_t(t)=-\tau^{-1}g^{\tau}(t)$, its time-derivative is given by
\begin{align*}
\frac{\mathrm{d}}{\mathrm{d}t}\mathbb{E}_1[\hat{u}]&=2\mathrm{Re}\left[\big(\tau\hat{u}_{ttt}+(\delta+\tau)|\xi|^2\hat{u}_t+|\xi|^2\hat{u}\big)\bar{\hat{u}}_{tt}\right]-\frac{|\xi|^2}{\tau^2}\int_0^tg^{\tau}(t-\eta)|\hat{u}(t,\xi)-\hat{u}(\eta,\xi)|^2\mathrm{d}\eta\\
&\quad+2|\xi|^2|\hat{u}_t|^2+\frac{2|\xi|^2}{\tau}\int_0^tg^{\tau}(t-\eta)\mathrm{Re}\left[\big(\hat{u}(t,\xi)-\hat{u}(\eta,\xi)\big)\bar{\hat{u}}_t(t,\xi)\right]\mathrm{d}\eta-\frac{g^{\tau}(t)}{\tau}|\xi|^2|\hat{u}|^2\\
&\quad+2|\xi|^2\big(g^{\tau}(t)-m\big)\mathrm{Re}(\hat{u}\bar{\hat{u}}_t)+\frac{2|\xi|^2}{\tau}\int_0^tg^{\tau}(t-\eta)\mathrm{Re}\big(\hat{u}(\eta,\xi)\bar{\hat{u}}_t(t,\xi)\big)\mathrm{d}\eta\notag\\
&\quad-2|\xi|^2\int_0^tg^{\tau}(t-\eta)\mathrm{Re}\big(\hat{u}(\eta,\xi)\bar{\hat{u}}_{tt}(t,\xi)\big)\mathrm{d}\eta.
\end{align*}
By using the equation in the inhomogeneous Cauchy problem \eqref{Inhomo-MGT-Memory}, one notices
\begin{align*}
2\mathrm{Re}\left[\big(\tau\hat{u}_{ttt}+(\delta+\tau)|\xi|^2\hat{u}_t+|\xi|^2\hat{u}\big)\bar{\hat{u}}_{tt}\right]=2\mathrm{Re}\big((-\hat{u}_{tt}+|\xi|^2g^{\tau}\ast\hat{u})\bar{\hat{u}}_{tt}\big)+2\delta|\xi|^2\mathrm{Re}(\widehat{\psi}_t^{\delta=0}\bar{\hat{u}}_{tt}).
\end{align*}
Then, associating with the fact
\begin{align*}
\frac{2|\xi|^2}{\tau}\int_0^tg^{\tau}(t-\eta)\mathrm{d}\eta=2|\xi|^2\big(m-g^{\tau}(t)\big),
\end{align*}
and observing some cancellations, we may derive
\begin{align}\label{Est-E1}
\frac{\mathrm{d}}{\mathrm{d}t}\mathbb{E}_1[\hat{u}]&=-2|\hat{u}_{tt}|^2+2|\xi|^2|\hat{u}_t|^2-\frac{|\xi|^2}{\tau^2}\int_0^tg^{\tau}(t-\eta)|\hat{u}(t,\xi)-\hat{u}(\eta,\xi)|^2\mathrm{d}\eta\notag\\
&\quad-\frac{g^{\tau}(t)}{\tau}|\xi|^2|\hat{u}|^2+2\delta|\xi|^2\mathrm{Re}(\widehat{\psi}_t^{\delta=0}\bar{\hat{u}}_{tt}).
\end{align}

For another, we introduce the second energy in the Fourier space
\begin{align*}
	\mb{E}_2[\hat{u}]&:=|\xi|^2|\hat{u}|^2+|\hat{u}_t|^2+2\tau\mathrm{Re}(\hat{u}_{tt}\bar{\hat{u}}_t)+|\xi|^2\int_0^tg^{\tau}(t-\eta)|\hat{u}(t,\xi)-\hat{u}(\eta,\xi)|^2\mathrm{d}\eta\\
	&\ \quad+\tau\big(g^{\tau}(t)-m\big)|\xi|^2|\hat{u}|^2,
\end{align*}
which shows via the analogous ways as those in the above that
\begin{align}\label{Est-E2}
	\frac{\mathrm{d}}{\mathrm{d}t}\mb{E}_2[\hat{u}]&=2\mathrm{Re}\big((|\xi|^2\hat{u}+\hat{u}_{tt}+\tau\hat{u}_{ttt})\bar{\hat{u}}_{t}\big)+2\tau|\hat{u}_{tt}|^2-\frac{|\xi|^2}{\tau}\int_0^tg^{\tau}(t-\eta)|\hat{u}(t,\xi)-\hat{u}(\eta,\xi)|^2\mathrm{d}\eta\notag\\
&\quad+2|\xi|^2\int_0^tg^{\tau}(t-\eta)\mathrm{Re}\left[\big(\hat{u}(t,\xi)-\hat{u}(\eta,\xi)\big)\bar{\hat{u}}_t(t,\xi)\right]\mathrm{d}\eta-g^{\tau}(t)|\xi|^2|\hat{u}|^2\notag\\
&\quad+2\tau \big(g^{\tau}(t)-m\big)|\xi|^2\mathrm{Re}(\hat{u}\bar{\hat{u}}_t)\notag\\
&=-2(\delta+\tau)|\xi|^2|\hat{u}_t|^2+2\tau|\hat{u}_{tt}|^2-g^{\tau}(t)|\xi|^2|\hat{u}|^2+2\delta|\xi|^2\mathrm{Re}(\widehat{\psi}_t^{\delta=0}\bar{\hat{u}}_{t})\notag\\
&\quad-\frac{|\xi|^2}{\tau}\int_0^tg^{\tau}(t-\eta)|\hat{u}(t,\xi)-\hat{u}(\eta,\xi)|^2\mathrm{d}\eta.
\end{align}
Then, we define a total energy $\mb{E}[\hat{u}]:=\mb{E}_1[\hat{u}]+k_0\mb{E}_2[\hat{u}]$ with a positive parameter $k_0$ to be chosen later. This total energy can be re-organized as follows:
\begin{align*}
	\mb{E}[\hat{u}]
	&=\frac{1}{k_0}\left[1+\tau\big(g^{\tau}(t)-m\big)\right]|\xi|^2|\hat{u}_t+k_0\hat{u}|^2+\left((\delta+\tau)-\frac{1}{k_0}\right)|\xi|^2|\hat{u}_t|^2+\tau|\hat{u}_{tt}+k_0\hat{u}_t|^2\\
	&\quad+k_0(1-k_0\tau)|\hat{u}_t|^2+\frac{|\xi|^2}{\tau}\int_0^tg^{\tau}(t-\eta)|\hat{u}(t,\xi)-\hat{u}(\eta,\xi)|^2\mathrm{d}\eta+g^{\tau}(t)|\xi|^2|\hat{u}|^2\\
	&\quad+\frac{\tau}{k_0}\big(m-g^{\tau}(t)\big)|\xi|^2|\hat{u}_t|^2+k_0|\xi|^2\int_0^tg^{\tau}(t-\eta)|\hat{u}(t,\xi)-\hat{u}(\eta,\xi)|^2\mathrm{d}\eta\\
	&\quad+2|\xi|^2\int_0^tg^{\tau}(t-\eta)\mathrm{Re}\left[\big(\hat{u}(t,\xi)-\hat{u}(\eta,\xi)\big)\bar{\hat{u}}_t(t,\xi)\right]\mathrm{d}\eta.
\end{align*}
	In the above, we used an important observation in the treatment of three terms in $\mb{E}[\hat{u}]$, precisely,
\begin{align*}
	&\underbrace{2|\xi|^2\mathrm{Re}(\hat{u}\bar{\hat{u}}_t)}_{\text{from }\mb{E}_1[\hat{u}]}+\underbrace{k_0\left[1+\tau\big(g^{\tau}(t)-m\big)\right]|\xi|^2|\hat{u}|^2}_{\text{from } \mb{E}_2[\hat{u}]}\\
	&=\frac{1}{k_0}\left[1+\tau\big(g^{\tau}(t)-m\big)\right]|\xi|^2|\hat{u}_t+k_0\hat{u}|^2-\frac{1}{k_0}|\xi|^2|\hat{u}_t|^2+\frac{\tau}{k_0}\big(m-g^{\tau}(t)\big)|\xi|^2|\hat{u}_t|^2\\
	&\quad+2|\xi|^2\int_0^tg^{\tau}(t-\eta)\mathrm{Re}\big(\hat{u}(t,\xi)\bar{\hat{u}}_t(t,\xi)\big)\mathrm{d}\eta.
\end{align*}
Clearly, Cauchy's inequality implies
\begin{align*}
	&2|\xi|^2\int_0^tg^{\tau}(t-\eta)\mathrm{Re}\left[\big(\hat{u}(t,\xi)-\hat{u}(\eta,\xi)\big)\bar{\hat{u}}_t(t,\xi)\right]\mathrm{d}\eta\\
	&\geqslant -k_0|\xi|^2\int_0^tg^{\tau}(t-\eta)|\hat{u}(t,\xi)-\hat{u}(\eta,\xi)|^2\mathrm{d}\eta-\frac{\tau}{k_0}\big(m-g^{\tau}(t)\big)|\xi|^2|\hat{u}_t|^2,
\end{align*}
 which leads to
 \begin{align}\label{Lower-E}
 \mb{E}[\hat{u}]
 &\geqslant \frac{1}{k_0}\left[1+\tau\big(g^{\tau}(t)-m\big)\right]|\xi|^2|\hat{u}_t+k_0\hat{u}|^2+\left((\delta+\tau)-\frac{1}{k_0}\right)|\xi|^2|\hat{u}_t|^2+\tau|\hat{u}_{tt}+k_0\hat{u}_t|^2\notag\\
 &\quad+k_0(1-k_0\tau)|\hat{u}_t|^2+\frac{|\xi|^2}{\tau}\int_0^tg^{\tau}(t-\eta)|\hat{u}(t,\xi)-\hat{u}(\eta,\xi)|^2\mathrm{d}\eta+g^{\tau}(t)|\xi|^2|\hat{u}|^2.
 \end{align}

From \eqref{Est-E1} as well as \eqref{Est-E2}, we may estimate
\begin{align*}
	\frac{\mathrm{d}}{\mathrm{d}t}\mb{E}[\hat{u}]&\leqslant 2(\tau k_0-1)|\hat{u}_{tt}|^2+2\big(1-(\delta+\tau)k_0\big)|\xi|^2|\hat{u}_t|^2+2\delta|\xi|^2\mathrm{Re}\big(\widehat{\psi}_t^{\delta=0}(\bar{\hat{u}}_{tt}+k_0\bar{\hat{u}}_t)\big)\\
	&\leqslant\left(\frac{\delta^2}{2-2\tau k_0}|\xi|^4+\frac{\delta^2k_0^2}{2(\delta+\tau)k_0-2}|\xi|^2\right)|\widehat{\psi}_t^{\delta=0}|^2\\
	&\leqslant C\delta|\xi|^2\langle\xi\rangle^2|\widehat{\psi}_t^{\delta=0}|^2,
\end{align*}
where the positive constant $C$ is independent of $\delta$, by fixing
\begin{align*}
k_0=\frac{\delta+2\tau}{2\tau(\delta+\tau)}\in\left(\frac{1}{\delta+\tau},\frac{1}{\tau}\right).
\end{align*}
 When $\xi\in\ml{Z}_{\intt}(\varepsilon_0)\cup\ml{Z}_{\extt}(N_0)$, because of the pairwise distinct characteristic roots, according to Section \ref{Section-MGT-Large-Time}, the time-derivative of the inviscid solution is represented by
\begin{align*}
	\widehat{\psi}^{\delta=0}_t&=\frac{\det(\mb{V}_1)}{\det(\mb{V})}\lambda_1\mathrm{e}^{\lambda_1t}+\frac{\det(\mb{V}_2)}{\det(\mb{V})}\lambda_2\mathrm{e}^{\lambda_2t}+\big(\cos(\lambda_{\mathrm{I}}t)\lambda_{\mathrm{R}}-\sin(\lambda_{\mathrm{I}}t)\lambda_{\mathrm{I}}\big)\frac{\det(\mb{V}_3)+\det(\mb{V}_4)}{\det(\mb{V})}\mathrm{e}^{\lambda_{\mathrm{R}}t}\\
	&\quad+i\big(\cos(\lambda_{\mathrm{I}}t)\lambda_{\mathrm{I}}+\sin(\lambda_{\mathrm{I}}t)\lambda_{\mathrm{R}}\big)\frac{\det(\mb{V}_3)-\det(\mb{V}_4)}{\det(\mb{V})}\mathrm{e}^{\lambda_{\mathrm{R}}t}.
\end{align*}
From asymptotic behavior derived in Section \ref{Section-MGT-Large-Time}, we obtain
\begin{align*}
\chi_{\intt}(\xi)|\widehat{\psi}^{\delta=0}_t|&\lesssim\chi_{\intt}(\xi)\mathrm{e}^{-c|\xi|^2t}\left(|\widehat{\psi}_0|+|\widehat{\psi}_1|+|\widehat{\psi}_2|\right),\\
\big(1-\chi_{\intt}(\xi)\big)|\widehat{\psi}^{\delta=0}_t|&\lesssim\big(1-\chi_{\intt}(\xi)\big)\mathrm{e}^{-c|\xi|^{-2}t}\left(|\widehat{\psi}_0|+|\widehat{\psi}_1|+\langle\xi\rangle^{-1}|\widehat{\psi}_2|\right).
\end{align*}
As a consequence,
\begin{align*}
\mb{E}[\hat{u}]&\leqslant C\delta|\xi|^2\langle\xi\rangle^2\int_0^t|\widehat{\psi}^{\delta=0}_t(\eta,\xi)|^2\mathrm{d}\eta\\
&\leqslant C\delta\left[\chi_{\intt}(\xi)\left(|\widehat{\psi}_0|^2+|\widehat{\psi}_1|^2+|\widehat{\psi}_2|^2\right)+\big(1-\chi_{\intt}(\xi)\big)\langle\xi\rangle^6\left(|\widehat{\psi}_0|^2+|\widehat{\psi}_1|^2+\langle\xi\rangle^{-2}|\widehat{\psi}_2|^2\right)\right],
\end{align*}
which is our desire uniform estimate in the diffusivity of sound $\delta$ and time $t$, simultaneously. From \eqref{Lower-E}, it is easy to see
\begin{align*}
|\xi|^2|\hat{u}_t+k_0\hat{u}|^2+|\hat{u}_{tt}+k_0\hat{u}_t|^2\leqslant C\mb{E}[\hat{u}].
\end{align*} By using the Hausdorff-Young inequality, we derive the global (in time) inviscid limit for the energy terms as follows:
\begin{align*}
&\|\,|D|(u_t+k_0u)(t,\cdot)\|_{L^{\infty}}+\|(u_{tt}+k_0u_t)(t,\cdot)\|_{L^{\infty}}\\
&\leqslant C\|\,|\xi|(\hat{u}_t+k_0\hat{u})(t,\xi)\|_{L^1}+ C\|(\hat{u}_{tt}+k_0\hat{u}_t)(t,\xi)\|_{L^1}\\
&\leqslant C\sqrt{\delta}\left\|\chi_{\intt}(\xi)\left(|\widehat{\psi}_0|+|\widehat{\psi}_1|+|\widehat{\psi}_2|\right)\right\|_{L^1}+C\sqrt{\delta}\left\|\big(1-\chi_{\intt}(\xi)\big)\left(\langle\xi\rangle^3|\widehat{\psi}_0|+\langle\xi\rangle^3|\widehat{\psi}_1|+\langle\xi\rangle^2|\widehat{\psi}_2|\right)\right\|_{L^1}\\
&\leqslant C\sqrt{\delta}\left\|\left(\widehat{\psi}_0,\widehat{\psi}_1,\widehat{\psi}_2\right)\right\|_{(L^{\infty})^3}+C\sqrt{\delta}\left\|\big(1-\chi_{\intt}(\xi)\big)\langle\xi\rangle^{2-s_0}\right\|_{L^1}\left\|\left(\langle\xi\rangle^{s_0+1}\widehat{\psi}_0,\langle\xi\rangle^{s_0+1}\widehat{\psi}_1,\langle\xi\rangle^{s_0}\widehat{\psi}_2\right)\right\|_{(L^{\infty})^3}\\
&\leqslant C\sqrt{\delta}\|(\psi_0,\psi_1,\psi_2)\|_{(H^{s_0+1}_1)^2\times H^{s_0}_1},
\end{align*}
where we restricted $n+1-s_0<-1$ to get the integrability
\begin{align*}
\left\|\big(1-\chi_{\intt}(\xi)\big)\langle\xi\rangle^{2-s_0}\right\|_{L^1}<+\infty.
\end{align*}
 Eventually, by the same method as the above one, from \eqref{Lower-E} again to get 
 \begin{align*}
 	|\hat{u}|\leqslant\frac{C}{\sqrt{g^{\tau}(t)}}|\xi|^{-1}\sqrt{\mb{E}[\hat{u}]},
 \end{align*}
  we obtain the local (in time) inviscid limit for the acoustic velocity potential
\begin{align*}
	\|u(t,\cdot)\|_{L^{\infty}}\leqslant C\sqrt{\delta}\mathrm{e}^{\frac{t}{2\tau}}\|(\psi_0,\psi_1,\psi_2)\|_{(H^{s_1+1}_1\cap \dot{H}^{-1}_1)^2\times (H^{s_1}_1\cap \dot{H}^{-1}_1)},
\end{align*}
with the restriction $n-s_1<-1$, in which the additional $\dot{H}^{-1}_1$ regularity of initial data comes from the factor $|\xi|^{-1}$ for small frequencies, and the growth rate $\mathrm{e}^{\frac{t}{2\tau}}$ comes from the fact $\sqrt{g^{\tau}(t)}$. Thus, our proof is complete.

\section{Final remarks: An application on the nonlinear Jordan-MGT equation in viscous hereditary fluids}\label{Section-Final-Remarks}\setcounter{equation}{0}
$\ \ \ \ $To end this manuscript, we will state an application of our optimal growth/decay estimates derived in Section \ref{Section-MGT-Large-Time} on the corresponding nonlinear acoustic waves model in viscous hereditary fluids. As we mentioned shortly in the introduction, the Cauchy problem for the viscous Jordan-MGT equation with memory has been investigated by \cite{Nikolic-Said=2021,Nikolic-Said=2021-02} recently, where the memory kernel was assumed to be a generally exponential decay function. Precisely, the global (in time) well-posedness for the Kuznetsov type model, i.e. \eqref{JMGT-Memory} with the nonlinearity $\partial_t(\frac{B}{2A}|\varphi_t|^2+|\nabla\varphi|^2)$, was proved in \cite{Nikolic-Said=2021} for $n\geqslant 3$ by the history framework of Dafermos transform combined with some energy bounds. Simultaneously, the authors of \cite{Nikolic-Said=2021-02} not only demonstrated the global (in time) solvability for the Westervelt type model, i.e. \eqref{JMGT-Memory} with the nonlinearity $(1+\frac{B}{2A})\partial_t(|\varphi_t|^2)$, for $n=3$ via energy methods, but also provided some decay estimates for energies via Fourier methods. 

Nevertheless, to the best of author's knowledge, the refined global (in time) behavior for the viscous Jordan-MGT equation with memory of type I in the whole space $\mb{R}^n$ seems open, especially, the global (in time) well-posedness in lower dimensions $n=1,2$ (the difficulties come from non-summable singularities in the linearized problem and Sobolev embeddings), the optimal large time estimates (cf. the sharpness question proposed in \cite[Remark 2]{Nikolic-Said=2021-02}), asymptotic profiles for large time. Later, we will state some results and enlightenment for these questions under the assumption \eqref{Memory-Kernel} on the memory kernel $g^{\tau}(t)$. Remark that without such hypothesis \eqref{Memory-Kernel}, the above questions are still challenging, which is beyond the scope of our paper.

Before doing these, let us employ the representation \eqref{Rep-wide-psi} associated with asymptotic behavior derived in Section \ref{Section-MGT-Large-Time} to obtain
\begin{align*}
|\partial_t^k\widehat{\psi}^{\delta>0}|&\lesssim\chi_{\intt}(\xi)|\xi|^{k-1}\mathrm{e}^{-c|\xi|^2t}\left(|\widehat{\psi}_0|+|\widehat{\psi}_1|+|\widehat{\psi}_2|\right)\\
&\quad+\big(1-\chi_{\intt}(\xi)\big)|\xi|^{k-1}\mathrm{e}^{-ct}\left(|\widehat{\psi}_0|+|\widehat{\psi}_1|+|\xi|^{-1}|\widehat{\psi}_2|\right)
\end{align*}
with $k=1,2$ for the Cauchy problem \eqref{MGT-Memory} with $\delta>0$. Therefore, the viscous solution also fulfills the following $(L^2\cap L^1)-L^2$ type estimates:
\begin{align*}
\|\partial_t^k\psi^{\delta>0}(t,\cdot)\|_{\dot{H}^s}\lesssim (1+t)^{\frac{1}{2}-\frac{s+k}{2}-\frac{n}{4}}\|(\psi_0,\psi_1,\psi_2)\|_{(H^{s+k-1}\cap L^1)^2\times(H^{s+k-2}\cap L^1)}
\end{align*}
with $k=1,2$ and $s\geqslant0$. Note that the decay rates and the regularity for the last data of these estimates coincide with those in \cite{Chen-Takeda=2023}. All steps in deriving global (in time) behavior of solutions, including global (in time) existence, optimal estimates and optimal leading term of small data solutions in \cite[Sections 3 and 4]{Chen-Takeda=2023}, because their proofs only depend on
\begin{itemize}
	\item the optimal decay rates of $(L^2\cap L^1)-L^2$ type estimates to construct suitable time-weighted Sobolev spaces;
	\item the regularity of the third data to estimate nonlinear terms in the $L^1$, $L^2$ as well as $\dot{H}^s$ norms via some fractional tools in the harmonic analysis;
	\item the representation of nonlinearity to reformulate a decomposition for the term $\varphi_t\varphi_{tt}$.
\end{itemize}
Fortunately, these three factors are the same as ours, which leads to the next results. Due to the similarity, we omit the detail of the proof.

We may demonstrate global (in time) existence of small data solution to the following viscous Jordan-MGT equation with memory:
\begin{align}\label{JMGT-Cauchy}
\begin{cases}
\displaystyle{\tau\varphi^{\mathrm{K}}_{ttt}+\varphi^{\mathrm{K}}_{tt}-\Delta\varphi^{\mathrm{K}}-(\delta+\tau)\Delta\varphi^{\mathrm{K}}_t+g^{\tau}\ast\Delta \varphi^{\mathrm{K}}_t=\partial_t\left(\frac{B}{2A}|\varphi^{\mathrm{K}}_t|^2+|\nabla \varphi^{\mathrm{K}}|\right)},&x\in\mb{R}^n,\ t>0,\\
(\varphi^{\mathrm{K}},\varphi^{\mathrm{K}}_t,\varphi^{\mathrm{K}}_{tt})(0,x)=(\varphi^{\mathrm{K}}_0,\varphi^{\mathrm{K}}_1,\varphi^{\mathrm{K}}_2)(x),&x\in\mb{R}^n,
\end{cases}
\end{align}
with $\tau>0$ and $\delta>0$, where the memory kernel $g^{\tau}(t)$ fulfills the condition \eqref{Memory-Kernel}. Precisely, there is a constant $\epsilon>0$ such that for all
\begin{align*}
(\varphi^{\mathrm{K}}_0,\varphi^{\mathrm{K}}_1,\varphi^{\mathrm{K}}_2)\in\ml{M}_s:=(H^{s+2}\cap L^1)\times(H^{s+1}\cap L^1)\times (H^s\cap L^1)
\end{align*} 
with $\|(\varphi^{\mathrm{K}}_0,\varphi^{\mathrm{K}}_1,\varphi^{\mathrm{K}}_2)\|_{\ml{M}_s}\leqslant \epsilon$, there is a uniquely determined Sobolev solution
\begin{align*}
\varphi^{\mathrm{K}}\in\ml{C}([0,+\infty),H^{s+2})\cap \ml{C}^1([0,+\infty),H^{s+1})\cap \ml{C}^2([0,+\infty),H^{s})
\end{align*}
for $s>\max\{\frac{n}{2}-1,0\}$ with $n\geqslant 1$, to the viscous Jordan-MGT equation of Kuznetsov type \eqref{JMGT-Cauchy}. Furthermore, the following estimates hold:
\begin{align*}
\|\varphi^{\mathrm{K}}(t,\cdot)\|_{L^2}&\lesssim\ml{D}_n(1+t)\|(\varphi^{\mathrm{K}}_0,\varphi^{\mathrm{K}}_1,\varphi^{\mathrm{K}}_2)\|_{\ml{M}_s},\\
\|\partial_t^{\ell}\varphi^{\mathrm{K}}(t,\cdot)\|_{L^2}&\lesssim (1+t)^{\frac{1}{2}-\frac{\ell}{2}-\frac{n}{4}}\|(\varphi^{\mathrm{K}}_0,\varphi^{\mathrm{K}}_1,\varphi^{\mathrm{K}}_2)\|_{\ml{M}_s}\ \ \mbox{when}\ \ \ell=1,2,\\
\|\partial_t^{\ell}\varphi^{\mathrm{K}}(t,\cdot)\|_{\dot{H}^{s+2-\ell}}&\lesssim (1+t)^{-\frac{1}{2}-\frac{s}{2}-\frac{n}{4}}\|(\varphi^{\mathrm{K}}_0,\varphi^{\mathrm{K}}_1,\varphi^{\mathrm{K}}_2)\|_{\ml{M}_s}\ \ \mbox{when}\ \ \ell=0,1,2.
\end{align*}
Its proof is based on the fixed point theory associated with the operator 
\begin{align*}
N^{\mathrm{K}}:\ \varphi^{\mathrm{K}}\to N^{\mathrm{K}}\varphi^{\mathrm{K}}:=\varphi^{\lin}+\int_0^tK_2(t-\eta,|D|)\partial_t\left(\frac{B}{2A}|\varphi^{\mathrm{K}}_t(\eta,x)|^2+|\nabla \varphi^{\mathrm{K}}(\eta,x)|\right)\mathrm{d}\eta
\end{align*}
with the kernel $K_2(t,|D|)$ for the third data, where $\varphi^{\lin}=\psi^{\delta>0}$ is the viscous solution to the linearized Cauchy problem \eqref{MGT-Memory} with $\delta>0$. Furthermore, one may derive the optimality for the last mentioned estimates. We conjugate that under the condition 
\begin{align*}
	P_{\varphi^{\mathrm{K}}_1+\tau\varphi^{\mathrm{K}}_2}-\tau\int_{\mb{R}^n}\left(\frac{B}{2A}|\varphi_1^{\mathrm{K}}(x)|^2+|\nabla\varphi^{\mathrm{K}}_0(x)|^2\right)\mathrm{d}x\neq0 
\end{align*}
the last estimates for $\varphi^{\mathrm{K}}$ are optimal for large time $t\gg1$. Moreover, the optimal leading term for the nonlinear problem \eqref{JMGT-Cauchy} also can be obtained by the same manner as the one in \cite{Chen-Takeda=2023}. Remark that this strategy also can be applied in the viscous Jordan-MGT equation of Westervelt type with memory of type I as follows:
\begin{align}\label{JMGT-Cauchy-W}
	\begin{cases}
		\displaystyle{\tau\varphi^{\mathrm{W}}_{ttt}+\varphi^{\mathrm{W}}_{tt}-\Delta\varphi^{\mathrm{W}}-(\delta+\tau)\Delta\varphi^{\mathrm{W}}_t+g^{\tau}\ast\Delta \varphi^{\mathrm{W}}_t=\left(1+\frac{B}{2A}\right)\partial_t(|\varphi^{\mathrm{W}}_t|^2)},&x\in\mb{R}^n,\ t>0,\\
		(\varphi^{\mathrm{W}},\varphi^{\mathrm{W}}_t,\varphi^{\mathrm{W}}_{tt})(0,x)=(\varphi^{\mathrm{W}}_0,\varphi^{\mathrm{W}}_1,\varphi^{\mathrm{W}}_2)(x),&x\in\mb{R}^n,
	\end{cases}
\end{align}
with $\tau>0$ and $\delta>0$, where the memory kernel $g^{\tau}(t)$ is given by \eqref{Memory-Kernel}. The proof is easier due to the lack of quadratic gradient nonlinearity $\partial_t(|\nabla\varphi|^2)$.

\appendix
\section{Optimal growth estimates for the inviscid MGT equation}\label{Appendix-A}
$\ \ \ \ $This appendix contributes to the optimal upper bound and lower bound estimates of the solution $\psi^{m=0,\delta=0}=\psi^{m=0,\delta=0}(t,x)$ to the inviscid MGT equation \eqref{MGT} with $\delta=0$ for large time in lower dimensions $n=1,2$. This result is a supplement to the previous literature \cite{Chen-Palmieri=2020,Kaltenbacher-Niko-2021} for the acoustic waves in inviscid fluids. It also may be regarded as the generalized results on the free wave equation ($m=\delta=\tau=0$) in the recent work \cite{Ikehata=2023}. 
\begin{theorem}\label{Prop-Inviscid-Optimal-Growth}
Suppose that initial data $(\psi_0,\psi_1,\psi_2)\in L^2\times (L^2\cap L^1)^2$. Then, the solution $\psi^{m=0,\delta=0}$ to the inviscid MGT equation \eqref{MGT} with $\delta=0$ in lower dimensions satisfies the following optimal growth estimates:
\begin{align*}
	\ml{D}_n(t)|P_{\psi_1+\tau\psi_2}|\lesssim\|\psi^{m=0,\delta=0}(t,\cdot)\|_{L^2}\lesssim\ml{D}_n(t)\|(\psi_0,\psi_1,\psi_2)\|_{L^2\times (L^2\cap L^1)^2}
\end{align*}
for $n=1,2$ and large time $t\gg1$, provided that $P_{\psi_1+\tau\psi_2}\neq0$.
\end{theorem}
\begin{proof}Our demonstration is fundamental and strongly motivated by the recent work \cite{Ikehata=2023}.
The trick is to re-organize the representation \eqref{Rep-Inviscid} as follows:
\begin{align*}
\widehat{\psi}^{m=0,\delta=0}&=\left(\frac{\cos(|\xi|t)}{1+\tau^2|\xi|^2}+\frac{\tau|\xi|\sin(|\xi|t)}{1+\tau^2|\xi|^2}+\frac{\tau^2|\xi|^2}{2(1+\tau^2|\xi|^2)}\mathrm{e}^{-\frac{t}{\tau}}\right)\widehat{\psi}_0\\
&\quad+\left(-\frac{\tau^3|\xi|\sin(|\xi|t)}{1+\tau^2|\xi|^2}-\frac{\tau^2\cos(|\xi|t)}{1+\tau^2|\xi|^2}+\frac{\tau^2}{2(1+\tau^2|\xi|^2)}\mathrm{e}^{-\frac{t}{\tau}}\right)\widehat{\psi}_2\\
&\quad+\frac{\sin(|\xi|t)}{|\xi|}\left(\widehat{\psi}_1+\tau\widehat{\psi}_2\right).
\end{align*}
On the right-hand side of the above formula, the first and second terms are the remainders producing bounded estimates only.
Hence, separating it into several parts, we arrive at
\begin{align*}
\|\psi^{m=0,\delta=0}(t,\cdot)\|_{L^2}&\lesssim\left\|\chi_{\intt}(\xi)\left(1+|\xi|^2\mathrm{e}^{-\frac{t}{\tau}}\right)\widehat{\psi}_0\right\|_{L^2}+\left\|\big(1-\chi_{\intt}(\xi)\big)\left(|\xi|^{-1}+\mathrm{e}^{-\frac{t}{\tau}}\right)\widehat{\psi}_0\right\|_{L^2}\\
&\quad+\left\|\chi_{\intt}(\xi)\left(1+\mathrm{e}^{-\frac{t}{\tau}}\right)\widehat{\psi}_2\right\|_{L^2}+\left\|\big(1-\chi_{\intt}(\xi)\big)\left(|\xi|^{-1}+|\xi|^{-2}\mathrm{e}^{-\frac{t}{\tau}}\right)\widehat{\psi}_2\right\|_{L^2}\\
&\quad+\left\|\frac{\sin(|\xi|t)}{|\xi|}\right\|_{L^2}\|(\psi_1,\psi_2)\|_{(L^1)^2}\\
&\lesssim\|(\psi_0,\psi_2)\|_{(L^2)^2}+\ml{D}_n(t)\|(\psi_1,\psi_2)\|_{(L^1)^2}
\end{align*}
for large time $t\gg1$, where we used the derived estimates in \cite[Equations (3.5) and (3.18)]{Ikehata=2023}. Similarly, by applying \cite[Proposition A.1]{Chen-Takeda=2023} and the triangle inequality $2|f+h|^2\geqslant |f|^2-2|h|^2$, one may deduce
\begin{align*}
\|\psi^{m=0,\delta=0}(t,\cdot)\|_{L^2}&\gtrsim\left\|\ml{F}^{-1}_{\xi\to x}\left(\frac{\sin(|\xi|t)}{|\xi|}\left(\widehat{\psi}_1+\tau\widehat{\psi}_2\right)\right)\right\|_{L^2}-\|(\psi_0,\psi_2)\|_{(L^2)^2}\\
&\gtrsim\ml{D}_n(t)|P_{\psi_1+\tau\psi_2}|-\|(\psi_0,\psi_2)\|_{(L^2)^2}
\end{align*}
for large time $t\gg1$. Then, recalling the growth property of $\ml{D}_n(t)$, our proof is complete.
\end{proof}

\section*{Acknowledgments}
 Wenhui Chen is supported in part by the National Natural Science Foundation of China (grant No. 12301270, grant No. 12171317), 2024 Basic and Applied Basic Research Topic--Young Doctor Set Sail Project (grant No. 2024A04J0016), Guangdong Basic and Applied Basic Research Foundation (grant No. 2023A1515012044).

\end{document}